\documentclass{article}
\usepackage{amsfonts}
\usepackage{amsmath,amssymb,amsmath}
\usepackage[dvips]{graphicx}
\usepackage{algorithm}
\usepackage{algorithmicx}
\usepackage{algpseudocode}
\usepackage{yhmath}
\usepackage{qtree}
\usepackage{xcolor}
\usepackage{epsfig}
\usepackage{lineno}
\usepackage{mathtools}
\usepackage{amsmath,bm}
\usepackage{amsthm}
\usepackage[T1]{fontenc}
\usepackage[utf8]{inputenc}
\usepackage{hyperref}
\usepackage{booktabs}
\usepackage{authblk}
\usepackage[toc]{appendix}
\usepackage{url}
\numberwithin{equation}{section}

\usepackage{verbatim}
\usepackage{geometry}

\usepackage[many]{tcolorbox}
\usetikzlibrary{shadows}
\usepackage{caption}
\newtcolorbox{shadedbox}{
  breakable,
  enhanced jigsaw,
  colback=white,
}

\newcommand{\be}{\begin{equation}}
\newcommand{\ee}{\end{equation}}
\newcommand{\ba}{\begin{array}}
\newcommand{\ea}{\end{array}}
\newcommand{\bea}{\begin{eqnarray*}}
\newcommand{\eea}{\end{eqnarray*}}
\newcommand{\bean}{\begin{eqnarray}}
\newcommand{\eean}{\end{eqnarray}}

\newtheorem{remark}{Remark}[section]

\makeatletter
\def\BState{\State\hskip-\ALG@thistlm}
\makeatother

\newcommand{\lc}{\mathrel{\raise2pt\hbox{${\mathop<\limits_{\raise1pt\hbox{\mbox{$\sim$}}}}$}}}
\newcommand{\gc}{\mathrel{\raise2pt\hbox{${\mathop>\limits_{\raise1pt\hbox{\mbox{$\sim$}}}}$}}}
\newcommand{\ec}{\mathrel{\raise1pt\hbox{${\mathop=\limits_{\raise2pt\hbox{\mbox{$\sim$}}}}$}}}

\newcommand{\blue}{\color{black}}
\newcommand{\red}{\color{black}}
\newcommand{\black}{\color{black}}
\newcommand{\orange}{\color{black}}

\begin{document}

\title{A Domain Decomposition Method for the Poisson-Boltzmann Solvation Models}

\author[1]{Chaoyu Quan}
\author[2,3]{Benjamin Stamm}
\author[4,5]{Yvon Maday \thanks{Accepted by SIAM Journal on Scientific Computing, December 5, 2018. Benjamin Stamm acknowledges the funding from the German Academic Exchange Service (DAAD) from funds of the “Bundesministeriums f{\"u}r Bildung und Forschung” (BMBF) for the project Aa-Par-T (Project-ID 57317909).
Yvon Maday and Chaoyu Quan acknowledge the funding from the PICS-CNRS and the PHC PROCOPE 2017 (Project No.37855ZK).
Chaoyu Quan acknowledges the financial support of the Fondation Sciences Math\'ematiques de Paris.
}}
\date{}
\affil[1]{\small Sorbonne Universit\'es, UPMC Univ Paris 06, UMR 7598, Laboratoire Jacques-Louis Lions, and Institut des Sciences du Calcul et des Donn\'ees, F-75005, Paris, France (\href{mailto:quan@ann.jussieu.fr}{quan@ann.jussieu.fr})}
\affil[2]{\small Center for Computational Engineering Science, RWTH Aachen University, Aachen, Germany (\href{mailto:best@mathcces.rwth-aachen.de}{best@mathcces.rwth-aachen.de})}
\affil[3]{\small Computational Biomedicine, Institute for Advanced Simulation IAS-5 and Institute of Neuroscience and Medicine INM-9, Forschungszentrum J\"ulich, Germany}
\affil[4]{\small Sorbonne Universit\'es, UPMC Univ Paris 06, UMR 7598, Laboratoire Jacques-Louis Lions, and Institut Universitaire de France, 75005, Paris, France (\href{mailto:maday@ann.jussieu.fr}{maday@ann.jussieu.fr})}
\affil[5]{\small Division of Applied Mathematics, Brown University, 182 George St, Providence, RI 02912, USA}

\maketitle

\begin{abstract}
In this paper, a domain decomposition method for the Poisson-Boltzmann (PB) solvation model that is widely used in computational chemistry is proposed.
This method, called ddLPB for short, solves the linear Poisson-Boltzmann (LPB) equation defined in $\mathbb R^3$ using the van der Waals cavity as the solute cavity.
The Schwarz domain decomposition method is used to formulate local problems  by decomposing the cavity into overlapping balls and only solving a set of coupled sub-equations in balls.
A series of numerical experiments is presented to test the robustness and the efficiency of this method including the comparisons with some existing methods.
We observe exponential convergence of the solvation energy with respect to the number of degrees of freedom which allows this method to reach the required level of accuracy  when coupling with quantum mechanical descriptions of the solute.
\end{abstract}


\section{Introduction}\label{sect:intro}
The properties of numerous charged bio-molecules and their complexes with other molecules are dependent on the dielectric permittivity and the ionic strength of their environment.
There are various methods to model ionic solution effects on molecular systems, which can be commonly divided into two broad categories
according to whether they employ an explicit or implicit solvation model. 
Explicit solvation models adopt molecular representations of both the solute and the solvent molecules, which produce accurate results, but are very expensive. 
Implicit solvation models adopt a microscopic treatment of the solute (with possibly a few solvent molecules), but characterize the solvent in terms of its macroscopic physical properties (for example, the solvent dielectric permittivity and the ionic strength).
This reduces greatly the computational cost compared to an explicit description of the solvent.
For this reason, implicit solvation models based on the Poisson-Boltzmann (PB) equation \cite{yoon1990boundary,nicholls1991rapid} are now widely-used, taking into account both the solvent (relative) dielectric permittivity and the ionic strength.
In this paper, we call these models the PB solvation models and we mention that the ESU-CGS (electrostatic units, centimetre-gram-second) system of units \cite{griffiths2008introduction} is used for all equations.

For the sake of simplicity, we consider the linear Poisson-Boltzmann (LPB) equation, which describes the electrostatic potential $\psi$ of the PB solvation model in the following form (see \cite{nicholls1991rapid})
\begin{equation}\label{eq:intro_linPB} 
- \nabla \cdot [\varepsilon(\mathbf x) \nabla \psi(\mathbf x)] + \bar\kappa(\mathbf x)^2 \psi(\mathbf x) = 4\pi\rho_{\rm M}(\mathbf x),\quad \mbox{in }\mathbb R^3,
\end{equation}
where $\varepsilon(\mathbf x)$ represents the space-dependent dielectric permittivity function, $\bar\kappa(\mathbf x)$ is the modified Debye-H{\"{u}}ckel parameter and $\rho_{\rm M}(\mathbf x)$ represents the known solute's charge distribution.
Usually, $\varepsilon(\mathbf x)$ has the following form
\begin{equation}
\begin{array}{r@{}l}
\varepsilon(\mathbf x) = \left\{
\begin{aligned}
&\varepsilon_{1} && \mbox{in $\Omega$},\\
&\varepsilon_{2} &&\mbox{in $\Omega^{\mathsf c} \coloneqq \mathbb R^3\backslash\overline \Omega$},
\end{aligned}
\right.
\end{array}\label{eq:eps}
\end{equation}
where $\varepsilon_1$ and $\varepsilon_2$ are respectively the solute dielectric permittivity and the solvent dielectric permittivity, $\Omega$ and $\Omega^{\mathsf c}$ represents respectively the solute cavity and the solvent region.
Furthermore, $\bar\kappa(\mathbf x)$ usually has the following form
\begin{equation}
\begin{array}{r@{}l}
\bar\kappa(\mathbf x) = \left\{
\begin{aligned}
&0 && \mbox{in $\Omega$},\\
&\sqrt{\varepsilon_2}\kappa &&\mbox{in $\Omega^{\mathsf c}$},
\end{aligned}
\right.
\end{array}\label{eq:kappa0}
\end{equation}
where $\kappa$ is the Debye-H\"uckel screening constant.
More details on the nonlinear Poisson-Boltzmann (NPB) equation and its linearization will be presented in Section \ref{sect:ionic_sol}.

Finally, we also mention two popular implicit solvation models as particular cases: the polarizable continuum model (PCM) \cite{cammi2007continuum,tomasi-2005,mennucci2010continuum} and the conductor-like screening model (COSMO) \cite{klamt1993cosmo}.
In the classical PCM, the solvent is represented as a polarizable continuous medium which is non-ionic, i.e., $\kappa = 0 $.
The COSMO is a reduced version of the PCM, where the solvent is represented as a conductor-like continuum. 
Both the PCM and the COSMO can be seen as two particular PB solvation models.
\black

\subsection{Previous work}
We recall three widely used methods for solving the LPB equation: the boundary element method (BEM), the finite difference method (FDM) and the finite element method (FEM), see \cite{lu2008recent} for a review.
As the names indicate, the BEM is based on solving an integral equation defined on the solute-solvent interface \cite{cances1998new}, {while the FDM and the FEM are implemented in some 3-dimensional big domain covering the solute molecule.}

In the BEM, the LPB equation is recast as some integral equations defined on the $2$-dimensional solute-solvent boundary \cite{yoon1990boundary,boschitsch2002fast,altman2009accurate,bajaj2011efficient}. 
To solve the integral equations, a surface mesh should be generated, for example, using the MSMS \cite{sanner-1996} or the NanoShaper \cite{zhang2006quality} etc.
The BEM is efficient to solve the LPB equation and some techniques can be used to accelerate the BEM solvers, including the fast multipole method \cite{zhang2015parallel} and the hierarchical ``treecode'' technique \cite{lu2008recent}.
\orange 
For instance, the PAFMPB solver \cite{lu2010afmpb, zhang2015parallel} developed by Lu et al. provides a fast calculation of the solvation energy, which uses the adaptive fast multipole method and achieves linear  complexity with respect to (w.r.t.) the number of mesh elements.
\black
Another interesting BEM solver, called TABI-PB \cite{geng2013treecode}, has been developed in the past several years, which uses the ``treecode'' technique.
However, the BEM has a limitation that it can not be easily generalized to solve the NPB equation.

To solve the general PB equation (linear or nonlinear), the FDM might be the most popular method.
Here, we list some successful FDM solvers: UHBD \cite{madura1995electrostatics}, DelPhi \cite{li2012delphi}, MIBPB by Wei's group \cite{chen2011mibpb} and APBS by Baker, Holst, McCammon et al. \cite{baker2001electrostatics,dolinsky2007pdb2pqr,jurrus2018improvements}.
In particular, the APBS is well-developed with many useful options and its popularity is still increasing.
In addition, there are some other contributions to the FDM for the PB equation \cite{nicholls1991rapid,fogolari2002poisson,qiao2006finite,li2012delphi, fisicaro2016generalized}.
In the FDM, a big box with grid is first taken, which covers the region of interest. 
Then, different types of boundary conditions can be chosen, such as zero, single Debye-H\"uckel, multiple Debye-H\"uckel and focusing boundary conditions (see the APBS documentation \cite{baker2001electrostatics,jurrus2018improvements}).
We mention that the cost of FDM can increase considerably with respect to the grid dimension, for example, when the grid dimension is $1000^{3}$ as mentioned in \cite{lu2008recent}.

Comparing to the BEM and the FDM, the FEM provides in general more flexibility for mesh refinement, more analysis of convergence and more selections of linear solvers \cite{lu2008recent}.
A rigorous solution and approximation theory of the FEM for the PB equation has been established in \cite{chen2007finite}.
Furthermore, the adaptive FEM developed by Holst et al. has tackled some of the most important issues of the PB equation \cite{holst2000adaptive,baker2000adaptive,chen2007finite, holst2012adaptive,Aksoylu2012}.
In addition, \blue the SDPBS and SMPBS web servers developed by Xie et al. for solving the size-modified PB equation have performed fast and efficiently \cite{xie2014new, ying2015new, jiang2015sdpbs, ying2018hybrid, xie2017smpbs}.
\black

In addition to the above methods, we mention the framework of particular domain decomposition methods for implicit solvation models (see also website \cite{ddfamily}).
In the past several years, a domain decomposition method for COSMO (called ddCOSMO) has been developed \cite{cances2013domain,lipparini2013fast,lipparini2014quantum,lipparini2014quantum2}.
This method is independent on mesh and grid, easy to implement, and about two orders of magnitude faster than the state of the art as demonstrated in \cite{lipparini2014quantum}.

The ddCOSMO method can be coupled with a quantum Hamiltonian  \cite{lipparini2014quantum,lipparini2014quantum2} or a polarizable force-field within molecular dynamics \cite{lipparini2015polarizable}.
Numerical tests of the method show linear scaling with respect to the number of atoms and
first results of these scaling properties of the ddCOSMO  in a simplified setting can be found in \cite{ciaramella2017analysis,ciaramella2018analysis}.
Recently, a similar discretization scheme for the classical PCM was proposed within the domain decomposition paradigm (called ddPCM) \cite{stamm2016new,gatto2017computation}.
Both the ddCOSMO and the ddPCM work for the solute cavity constituted by overlapping balls, such as the van der Waals (VDW) cavity and the solvent accessible surface (SAS) cavity \cite{lee-1971,lee-1977}.
In the case of the PCM based on the ``smooth'' molecular surface, i.e., based on the solvent excluded surface (SES) \cite{Connolly-1983, Quan2016760},  another domain decomposition method has been proposed in \cite{quan2018domain}, which is called the ddPCM-SES.

Inspired by the previous work mentioned above, we develop a particular domain decomposition method for the PB solvation model  (called ddLPB), to solve the LPB equation in $\mathbb R^3$. 

\subsection{ddLPB}
In fact, the LPB equation \eqref{eq:intro_linPB} consists of a Poisson equation defined in the bounded solute cavity $\Omega$ and a homogeneous screened Poisson (HSP) equation defined in the unbounded solvent region $\Omega^{\mathsf c}$, which are coupled by some jump conditions on the interface $\Gamma\coloneqq \partial \Omega$.

\begin{figure}
\centering
\includegraphics[width = 5in]{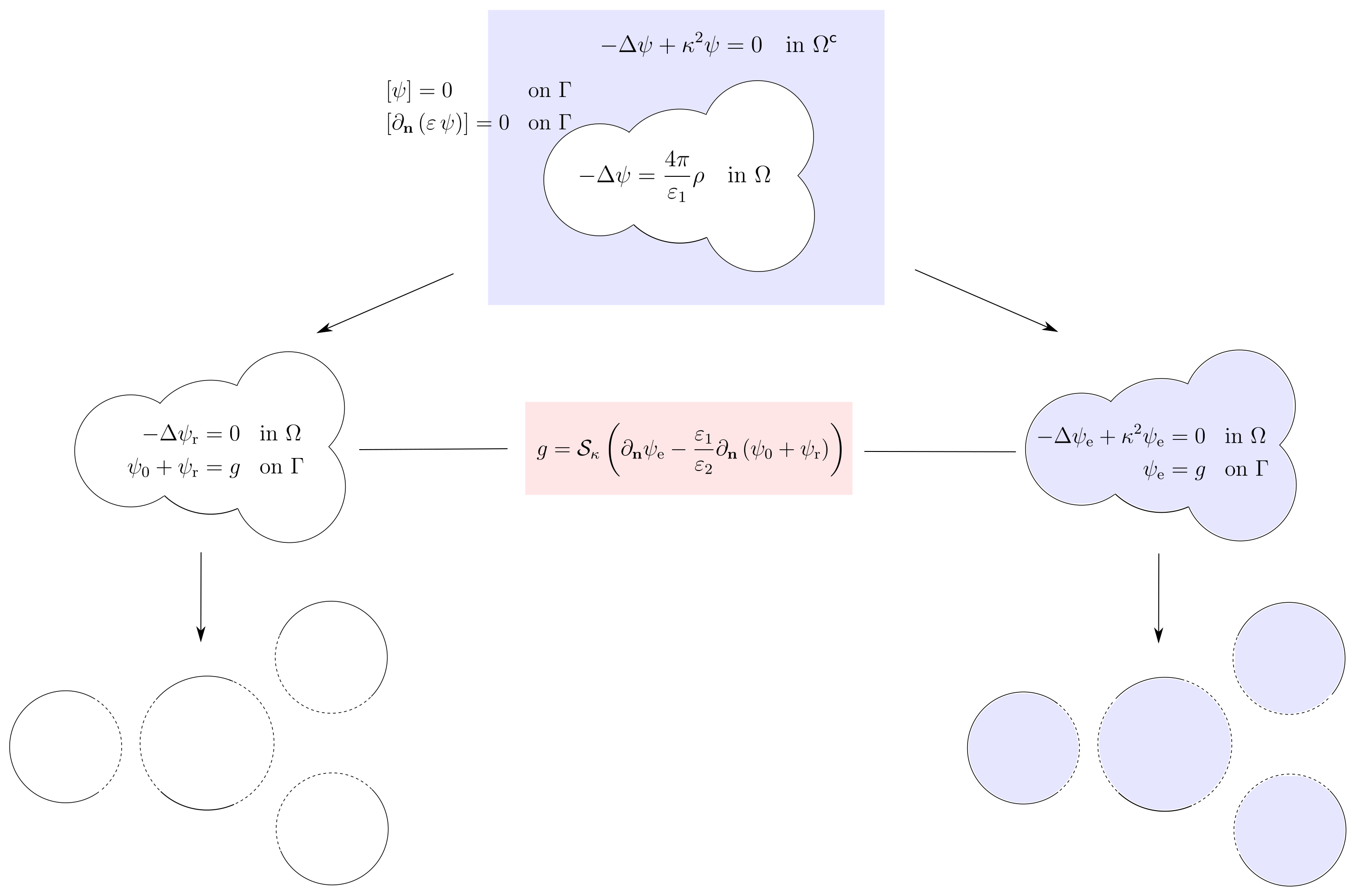}
\caption{Schematic diagram of the ddLPB.}\label{fig:scheme}
\end{figure}

To solve this problem, we first transform the Poisson equation in \eqref{eq:intro_linPB} into the following Laplace equation of $\psi_{\rm r}\coloneqq \psi - \psi_0$,
\begin{equation}\label{eq:contribution1}
-\Delta \psi_{\rm r} = 0,\quad \mbox{in } \Omega,
\end{equation}
where $\psi_{\rm r}$ is called the reaction potential and $\psi_0$ satisfies $-\Delta \psi_0 = \frac{4\pi}{\varepsilon_1} \rho_{\rm M}$ in $\mathbb R^3$.
Then, according to the potential theory, the electrostatic potential $\psi|_{\Omega^{\mathsf c}}$ can be represented as a single-layer potential (an exterior Dirichlet problem), which simultaneously gives an extended potential $\psi_{\rm e}$ satisfying the following HSP equation defined now in $\Omega$ (an interior Dirichlet problem)
\begin{equation}\label{eq:contribution2}
- \Delta\psi_{\rm e}(\mathbf x) + \kappa^2 \psi_{\rm e}(\mathbf x)= 0, \quad \mbox{in $\Omega$}.
\end{equation}
Based on the classical jump-conditions (see Section \ref{sect:tran}) of $\psi$ on the solute-solvent boundary, a coupling condition (see Figure \ref{fig:scheme}) between the Laplace equation \eqref{eq:contribution1} and the extended HSP equation \eqref{eq:contribution2} arises through an auxiliary function $g$ defined by
\begin{equation}\label{eq:contribution3}
g= \mathcal S_\kappa\left(\partial_{\mathbf n}\psi_{\rm e} - \frac{\varepsilon_1}{\varepsilon_2}\partial_{\mathbf n}\left(\psi_0+\psi_{\rm r} \right)\right),\quad \mbox{on $\Gamma$},
\end{equation}
where $\mathcal S_\kappa: H^{-\frac 1 2}(\Gamma) \rightarrow H^{\frac 1 2}(\Gamma)$ denotes the single-layer operator on $\Gamma$ ($\mathcal S_\kappa$ is defined in Section \ref{sect:tran}).
Here, $H^{-\frac 1 2}(\Gamma)$ and $H^{\frac 1 2}(\Gamma)$ denote the usual Sobolev spaces of order $\pm \frac 1 2$ on $\Gamma$, see  \cite{adams2003sobolev}. \black 
The initial problem defined in $\mathbb R^3$ is therefore transformed into two equations \eqref{eq:contribution1}--\eqref{eq:contribution2} coupled through $g$ in Eq. \eqref{eq:contribution3}.

Considering the fact that the solute cavity is commonly modeled as a union of overlapping balls, a particular Schwarz domain decomposition method (called ddLPB) can be used to solve Eqs \eqref{eq:contribution1}--\eqref{eq:contribution2} by respectively solving a group of coupled sub-equations in balls.
The main idea of this domain decomposition method is illustrated in Figure \ref{fig:scheme}.
Ultimately, only a Laplace solver and a HSP solver in the unit ball need to be developed for the local Laplace sub-equations and the local HSP sub-equations.
Each solver uses the spectral method for the corresponding PDE, where the spherical harmonics are taken as basis functions in the angular direction of the spherical coordinate system.

The ddLPB provides a new discretization of the LPB equation and has its own features.
In fact, this method is initially designed for quantum calculations, which usually require the accurate electrostatic solvation energy and the derivatives w.r.t. the atom positions. 
This method does not rely on mesh nor grid, but only on the Lebedev quadrature points \cite{lebedev1999quadrature} on $2$-dimensional spheres.
Therefore, it will be convenient to apply the ddLPB in molecular dynamics, without remeshing molecular surface as in the BEM.
The computation of forces becomes also very natural as the spheres are centered at the nuclear positions.
In the numerical tests, we will further show that the ddLPB is numerically robust and efficient.
\black

\subsection{Outline}
In Section \ref{sect:ionic_sol}, we introduce the derivation of the PB equation as well as its linearization.
Then, in Section \ref{sect:tran}, we transform the original LPB equation defined in $\mathbb R^3$ into the Laplace equation and the HSP equation both defined in the bounded solute cavity, as briefly outlined above.
We present a global strategy for solving the transformed problem.
This strategy involves solving the Laplace equation and the HSP equation defined in the solute cavity, which will be presented in Section \ref{sect:strategy} using a particular domain decomposition method.
In Section \ref{sect:solver}, we develop a Laplace solver and a HSP solver in the unit ball to solve the local equations.
After that, in Section \ref{sect:reform_disc}, we reformulate the coupling conditions and deduce a global linear system to be finally solved.
In Section \ref{sect:numerical}, we present some numerical results about the ddLPB.
In the last section, we draw some conclusions.

\section{PB solvation model} \label{sect:ionic_sol}
In this section, we introduce the well-known PB equation and its linearization, which describe the electrostatic potential in the implicit solvation model with ionic solutions. 

\blue
The space $\mathbb R^3$ is simply divided into the solute cavity and the solvent region, as introduced in Eqs \eqref{eq:intro_linPB} -- \eqref{eq:kappa0}.
Three types of molecular surfaces are mostly used to define the solute-solvent interface: the VDW surface, the SAS and the SES.
Both the VDW surface and the SAS are the boundary of the union of balls (respectively the VDW-balls and the SAS-balls), while the geometrical structure of SES is more complicated, see \cite{Quan2016760,quan2017meshing} for a thorough characterization.
In practice, the scaled VDW surface is often used, where each VDW-radius is multiplied by a scalar factor such as $1.1\sim1.2$ which is a common approach.
For the rest of this article we will limit the development to VDW-cavities.
Note that without any further difficulty, the ddLPB method also work for the scaled VDW-cavity and the SAS-cavity. 
\black


\subsection{Poisson-Boltzmann equation}\label{sect:PB}
In the PB solvation model, the solvent is represented by a polarizable and ionic continuum.
The freedom of the ions to move in the solution is accounted for by Boltzmann statistics. 
That is to say, the Boltzmann equation is used to calculate the local ion density $c_i$ of the $i$-th type of ion as follows
\begin{equation}\label{eq:boltzmann}
\displaystyle c_i = c^\infty_i  \exp\left({\frac{-W_i}{k_{\rm B}T}}\right),
\end{equation}
where
$c^\infty_i$ is the bulk ion concentration at an infinite distance from the solute molecule,
$W_i$ is the work required to move the $i$-th type of ion to a given position from an infinitely far distance,
$k_{\rm B}$ is the Boltzmann constant,
$T$ is the temperature in Kelvins (K).
The electrostatic potential $\psi$ of a general implicit solvation model is described originally by the Poisson equation as follows
\begin{equation}\label{eq:pde} 
-\nabla\cdot \varepsilon(\mathbf x) \nabla \psi(\mathbf x) = 4\pi \rho(\mathbf x), \quad \mbox{in } \mathbb R^3,
\end{equation}
where $\psi(\mathbf x) = O( \frac 1 {|\mathbf x|})$ as $|\mathbf x|\rightarrow \infty$.
Here, $\varepsilon(\mathbf x)$ represents the space-dependent dielectric permittivity and $\rho(\mathbf x)$ represents the charge distribution of the solvated system.
Given the solute's charge distribution $\rho_{\rm M}$ and the ionic distribution $c_i$ in \eqref{eq:boltzmann}, we can derive the PB equation from Eq. \eqref{eq:pde} as follows (see \cite{fogolari2002poisson})
\begin{equation}\label{eq:nonlinPB}
\displaystyle - \nabla \cdot [\varepsilon(\mathbf x) \nabla \psi(\mathbf x)] = 4\pi\rho_{\rm M}(\mathbf x) + \sum_i \,z_i \, e\, c^\infty_i  \, \exp\left(\frac{-z_i e\psi(\mathbf x)}{k_{\rm B}T}\right) \, \chi_{\Omega^{\mathsf c}}(\mathbf x),
\end{equation}
where $z_i e$ is the charge of the $i$-th type of ion, $e$ is the elementary charge and $\chi_{\Omega^{\mathsf c}}$ is the characteristic function of the solvent region $\Omega^{\mathsf c}$.

\begin{figure}
\centering
\includegraphics[width = 2.2 in, trim = {0.1in 0.1in 0.1in 0.1in},clip]{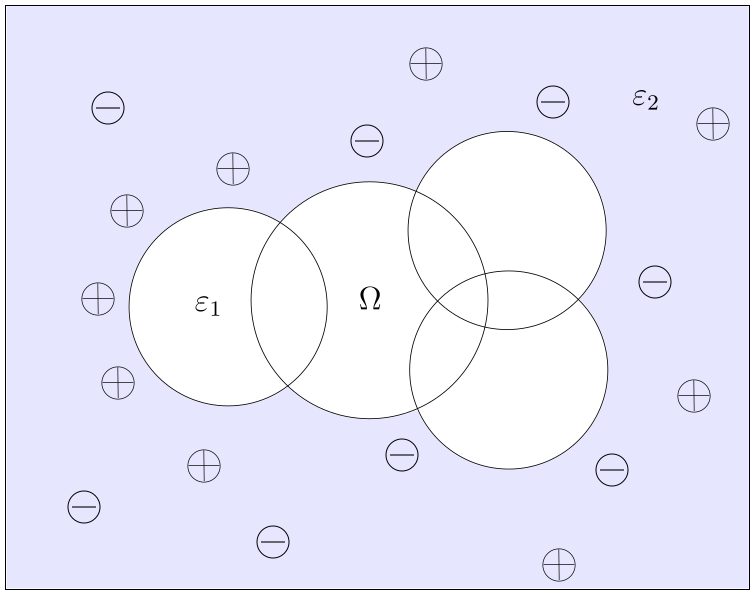}
\caption{2D schematic diagram of the implicit solvation model with ionic solutions, i.e., the PB solvation model.}\label{fig:pb}
\end{figure}

In the PB solvation model with a $1:1$ electrolyte, there are two types of ions respectively with charge $+e$ and $-e$ (see Figure \ref{fig:pb} for a schematic diagram).
With the assumption that $\psi$ satisfies the low potential condition, i.e., $\left | \frac{e\psi}{k_{\rm B}T} \right | \ll 1$, the NPB equation \eqref{eq:nonlinPB} can be linearized to (see \cite{nicholls1991rapid} for this form)
\begin{equation}\label{eq:linpb}
- \nabla \cdot [\varepsilon(\mathbf x) \nabla \psi(\mathbf x)] + \bar\kappa(\mathbf x)^2 \psi(\mathbf x) = 4\pi\rho_{\rm M}(\mathbf x),
\end{equation}
where $\psi$ is determined by the data $\varepsilon(\mathbf x)$, $\bar\kappa(\mathbf x)$ and $\rho_{\rm M}(\mathbf x)$ that are introduced in Section \ref{sect:intro}.

\begin{remark}
When the ionic solution has more than two types of ions, the nonlinear Poisson-Boltzmann equation can still be linearized to the form \eqref{eq:linpb}. 
In order to obtain a simpler expression of the modified Debye-H\"uckel parameter $\bar\kappa(\mathbf x)$, we consider the $1:1$ electrolyte in this paper.
\end{remark}

In the definition \eqref{eq:eps} of $\varepsilon(\mathbf x)$, the solute (relative) dielectric permittivity $\varepsilon_1$ should theoretically be set to $1$ as in the vacuum (for example, in \cite{cances1998new}). 
However, values different from $\varepsilon_1=1$ might be used.
For example, in \cite{nicholls1991rapid}, the authors claim to obtain better approximations with the empirical value $\varepsilon_1=2$.
The solvent dielectric permittivity $\varepsilon_2$ is determined by the solvent as well as the temperature, for example, $\varepsilon_2=78.54$ for water at the room temperature $25^\circ{\rm C}$.
The modified Debye-H\"uckel parameter in the implicit solvation model with a $1:1$ electrolyte is taken as
\begin{equation}
\begin{array}{r@{}l}
\bar\kappa(\mathbf x) = \left\{
\begin{aligned}
&0 && \mbox{in $\Omega$},\\
&\sqrt{\varepsilon_2}\kappa &&\mbox{in $\Omega^{\mathsf c}$},
\end{aligned}
\right.
\end{array}\label{eq:kappa}
\end{equation}
where $\kappa$ is the Debye-H\"uckel screening constant representing the attenuation of interactions due to the presence of ions in the solvent region, which is related to the ionic strength $I$ of the ionic solution according to (see \cite{nicholls1991rapid} and \cite[Section 1.4]{lamm2003poisson} for the following formula)
\begin{equation}\label{intro:eq:kappa}
\kappa^2 = \frac{8\pi e^2 N_{\rm A}I}{1000 \varepsilon_2 k_{\rm B} T},
\end{equation}
where $N_A$ is the Avogadro constant.

Furthermore, it is usually assumed that the solute's charge distribution $\rho_{\rm M}$ is supported in $\Omega$. 
For example, for a classical description of the solute, $\rho_{\rm M}$ is given by the sum of $M$ point charges in the following form
\begin{equation}\label{eq:rho}
\rho_{\rm M}(\mathbf x) = \sum_{i = 1}^M q_i\,\delta(\mathbf x - \mathbf x_i),
\end{equation}
where $M$ is the number of solute atoms, $q_i$ represents the (partial) charge carried on the $i$th atom with center $\mathbf x_i$, $\delta$ is the Dirac delta function.
For a quantum description of the solute, $\rho_{\rm M}$ consists of a sum of classical nuclear charges and the electron charge density.

\section{Problem transformation}\label{sect:tran}
In this section, we first introduce the integral representation of the LPB equation in the potential theory. 
Based on this, we then transform the original electrostatic problem to two coupled equations restricted to the (bounded) solute cavity.

\subsection{Problem setting}
The LPB equation can be divided into two equations: first, the Poisson equation in the solute cavity and second, the HSP equation in the solvent region. That is to say, the problem is recast in the following form
\begin{equation}\label{eq:pb_div}
\begin{array}{r@{}l}
\left\{
\begin{aligned}
&- \Delta \psi(\mathbf x)  = \frac{4\pi}{\varepsilon_1}\rho_{\rm M}(\mathbf x) && \mbox{in $\Omega$},\\
&- \Delta\psi(\mathbf x) + \kappa^2 \psi(\mathbf x) = 0 &&\mbox{in $\Omega^{\mathsf c}$},
\end{aligned}
\right.
\end{array}
\end{equation}
with two classical jump-conditions
\begin{equation}\label{eq:jc}
\begin{array}{r@{}l}
\left\{
\begin{aligned}
\left[\psi\right] &= 0 && \mbox{on $\Gamma$},\\
\left[\partial_{\mathbf n} \left(\varepsilon\,\psi\right)\right]&= 0 &&\mbox{on $\Gamma$},
\end{aligned}
\right.
\end{array}
\end{equation}
where $\Gamma \coloneqq \partial \Omega$ is the solute-solvent boundary, $\mathbf n$ is the unit normal vector on $\Gamma$ pointing outwards with respect to $\Omega$ and $\partial_{\mathbf n} = \mathbf n\cdot \nabla $ is the notation of normal derivative. $[\psi]$ represents the jump (inside minus outside) of the potential and $\left[\partial_{\mathbf n} \left(\varepsilon\,\psi\right)\right]$ represents the jump of the normal derivative of the electrostatic potential multiplied by the dielectric permittivity.

\subsection{Necessary tools from the potential theory}
The free-space Green's function of the operator $-\Delta$ is given as
\begin{equation}
G(\mathbf x,\mathbf y) = \frac{1}{4\pi |\mathbf x - \mathbf y|},\quad \forall \mathbf x,\mathbf y\in\mathbb R^3,
\end{equation}
and similarly, the free-space Green's function of the operator $-\Delta+\kappa^2$ is given as
\begin{equation}
G_\kappa(\mathbf x,\mathbf y) = \frac{\exp\left({-\kappa|\mathbf x - \mathbf y|}\right)}{4\pi |\mathbf x - \mathbf y|},\quad \forall \mathbf x,\mathbf y\in\mathbb R^3,
\end{equation}
which yields 
\begin{equation}
-\Delta_{\mathbf x} G(\mathbf x,\mathbf y) = \delta(\mathbf x-\mathbf y),\quad \forall \mathbf y \in \mathbb R^3,
\end{equation}
and
\begin{equation}
-\Delta_{\mathbf x} G_\kappa(\mathbf x,\mathbf y) + \kappa^2 G_\kappa(\mathbf x,\mathbf y) = \delta(\mathbf x-\mathbf y),\quad \forall \mathbf y \in \mathbb R^3.
\end{equation}

In the solute cavity $\Omega$, we define the reaction potential $\psi_{\rm r} \coloneqq \psi - \psi_0$, where $\psi_0$ is the potential generated by $\rho_{\rm M}$ in vacuum written as
\begin{equation}
\psi_0 =  \sum_{i = 1}^{M}\frac{q_i}{\varepsilon_1|\mathbf x - \mathbf x_i|},
\end{equation}
satisfying $-\Delta \psi_0 = \frac{4\pi}{\varepsilon_1} \rho_{\rm M}$ in $\mathbb R^3$.
Then, $\psi_{\rm r}$ is harmonic in $\Omega$, that is,
\begin{equation}\label{eq:psi_r}
- \Delta \psi_{\rm r} = 0,\quad \mbox{in } \Omega,
\end{equation}
which yields the following integral equation
\begin{equation} \label{eq:S_1}
\psi_{\rm r} (\mathbf x) = {\mathcal {\widetilde S}}\sigma_{\rm r} (\mathbf x)\coloneqq \int_\Gamma \frac{\sigma_{\rm r}(\mathbf y)}{4\pi |\mathbf x - \mathbf y|}  ,\quad \forall \mathbf x \in \Omega,
\end{equation}
where $\sigma_{\rm r}$ is some function in $H^{-\frac 1 2}(\Gamma)$ and $\mathcal {\widetilde S}: H^{-\frac 1 2}(\Gamma) \rightarrow H^{1}(\mathbb R^3\backslash \Gamma)$ is the single-layer potential associated with $G$.

Furthermore, according to the HSP equation in \eqref{eq:pb_div}, the electrostatic potential in the solvent region $\Omega^{\mathsf c}$ can be represented by
\begin{equation}\label{eq:S_2}
\psi|_{\Omega^{\mathsf c}} (\mathbf x) = {\mathcal {\widetilde S}}_\kappa \sigma_{\rm e} (\mathbf x) \coloneqq \int_\Gamma \frac{\exp\left({-\kappa |\mathbf x-\mathbf y|}\right)\sigma_{\rm e}(\mathbf y)}{4\pi |\mathbf x - \mathbf y|} ,\quad \forall \mathbf x \in \Omega^{\mathsf c}, 
\end{equation}
where $\sigma_{\rm e}$ is another function in $H^{-\frac 1 2}(\Gamma)$ and $\mathcal {\widetilde S}_\kappa: H^{-\frac 1 2}(\Gamma) \rightarrow H^{1}(\mathbb R^3\backslash \Gamma)$ is the single-layer potential associated with $G_\kappa$.
Here, we also introduce the single-layer operator $\mathcal S_\kappa: H^{-\frac 1 2}(\Gamma) \rightarrow H^{\frac 1 2}(\Gamma)$ defined by
\begin{equation}
\mathcal S_\kappa \sigma_{\rm e} (\mathbf x) \coloneqq \int_\Gamma \frac{\exp\left({-\kappa |\mathbf x-\mathbf y|}\right)\sigma_{\rm e}(\mathbf y)}{4\pi |\mathbf x - \mathbf y|},\quad \forall \mathbf x\in \Gamma,
\end{equation}
which is an invertible operator (this is true according to the proof of the invertibility of the single-layer operator for the Helmholtz equation, see \cite[Corollary 7.26]{L.Banjai} and \cite[Theorem 3.9.1]{Sauter2011}).
The invertibility of $\mathcal S_\kappa$ implies that $\sigma_{\rm e}$ can be characterized as $\sigma_{\rm e} = \mathcal S_\kappa^{-1} \psi|_\Gamma$.

\subsection{Transformation}
We will now transform the original problem defined in $\mathbb R^3$ equivalently to two coupled equations both defined in the solute cavity.

According to the continuity of the single-layer potential $\mathcal {\widetilde S}_\kappa$ across the interface \cite{Sauter2011}, we can artificially extend the electrostatic potential $\psi|_{\Omega^{\mathsf c}}$ from $\Omega^{\mathsf c}$ to $\Omega$ as follows
\begin{equation}\label{eq:S_sigma2}
\psi_{\rm e} (\mathbf x) \coloneqq {\mathcal {\widetilde S}}_\kappa \sigma_{\rm e} (\mathbf x) =\int_\Gamma \frac{\exp\left({-\kappa |\mathbf x-\mathbf y|}\right)\sigma_{\rm e}(\mathbf y)}{4\pi |\mathbf x - \mathbf y|}, \quad \forall \mathbf x\in \Omega,
\end{equation}
where $\psi_{\rm e}$ is called the extended potential in this paper.
As a consequence, $\psi_{\rm e}$ satisfies the same HSP equation as $\psi|_{\Omega^{\mathsf c}}$, but defined on $\Omega$, as follows
\begin{equation}\label{eq:psi_e}
-\Delta \psi_{\rm e} (\mathbf x) +\kappa^2 \psi_{\rm e} (\mathbf x)= 0,\quad \mbox{in~}\Omega,
\end{equation} 
with the same Dirichlet boundary conditions on $\Gamma$.
Furthermore, from \cite[Theorem 3.3.1]{Sauter2011}, we have a relation among $\sigma_{\rm e}$ and the normal derivatives of $\psi_{\rm e}$ and $\psi|_{\Omega^{\mathsf c}}$ on $\Gamma$ as follows
\begin{equation}\label{eq:relation}
\sigma_{\rm e} = \partial_{\mathbf n}\psi_{\rm e}|_{\Omega} - \partial_{\mathbf n}\psi|_{\Omega^{\mathsf c}},\quad \mbox{on }\Gamma.
\end{equation}

As introduced above, Eq. \eqref{eq:psi_r} of $\psi_{\rm r}$ and \eqref{eq:psi_e} of $\psi_{\rm e}$ are two PDEs defined on $\Omega$, that are derived from the original LPB equation \eqref{eq:pb_div}. 
As a consequence, it is sufficient to couple these two equations.
According to $[\psi] = 0$ on $\Gamma$ and the continuity of $\mathcal {\widetilde S}_\kappa$ across $\Gamma$ \cite{Sauter2011}, we then deduce a first coupling condition
\begin{equation}
\psi_0 + \psi_{\rm r}  = \psi_{\rm e},\quad \mbox{on }\Gamma.
\end{equation}
Further, combining Eq. \eqref{eq:relation} with the second equation of the jump conditions \eqref{eq:jc}, i.e.,
\begin{equation}
\varepsilon_1\partial_{\mathbf n}\psi|_{\Omega} - \varepsilon_2\partial_{\mathbf n}\psi|_{\Omega^{\mathsf c}} = 0,\quad \mbox{on }\Gamma,
\end{equation}
we deduce another coupling condition
\begin{equation}
\sigma_{\rm e}= \partial_{\mathbf n}\psi_{\rm e} - \frac{\varepsilon_1}{\varepsilon_2}\partial_{\mathbf n}\left(\psi_0+\psi_{\rm r} \right),\quad \mbox{on }\Gamma.
\end{equation}

In summary, the original problem \eqref{eq:pb_div} is transformed into the following two equations defined on $\Omega$
\begin{equation}\label{eq:tran_prob}
\begin{array}{r@{}l}
\left\{
\begin{aligned}
&- \Delta\psi_{\rm r}(\mathbf x) = 0 && \mbox{in $\Omega$},\\
&- \Delta\psi_{\rm e}(\mathbf x) + \kappa^2 \psi_{\rm e}(\mathbf x)= 0 &&\mbox{in $\Omega$},
\end{aligned}
\right.
\end{array}
\end{equation}
with two coupling conditions on $\Gamma$ given by
\begin{equation}\label{eq:tran_bd}
\begin{array}{r@{}l}
\left\{
\begin{aligned}
&\psi_0 + \psi_{\rm r}  = \psi_{\rm e} && \mbox{on }\Gamma,\\
&\sigma_{\rm e} = \partial_{\mathbf n}\psi_{\rm e} - \frac{\varepsilon_1}{\varepsilon_2}\partial_{\mathbf n}\left(\psi_0+\psi_{\rm r} \right) &&\mbox{on $\Gamma$},
\end{aligned}
\right.
\end{array}
\end{equation}
where $\sigma_{\rm e}$ is the charge density generating $\psi_{\rm e}$, as presented in \eqref{eq:S_sigma2}.
The second equation of \eqref{eq:tran_bd} is also equivalent to 
\begin{equation}\label{eq:tran_bd2}
\varepsilon_2\psi_{\rm e} 
+
\mathcal S_\kappa
\left(
\varepsilon_1\partial_{\mathbf n}\psi_{\rm r}
-
\varepsilon_2\partial_{\mathbf n}\psi_{\rm e}
\right)
=
-\varepsilon_1
\mathcal S_\kappa
\left(\partial_{\mathbf n}\psi_0\right),
\quad \mbox{on $\Gamma$},
\end{equation}
which is derived from letting $\mathcal S_\kappa$ act on both sides of the equation.

\begin{remark}
Eqs \eqref{eq:tran_prob} -- \eqref{eq:tran_bd2} are equivalent to the integral equation formulations (IEF) in  \cite{cances1997new} and \cite{geng2013treecode}.
The reason why we do the above transformation is that both the Laplace and the HSP equations in \eqref{eq:tran_prob} can be solved efficiently using a particular domain decomposition method,  see Section \ref{subsect:dd} and Section \ref{sect:solver}.
The main computational cost will be spent on the coupling conditions \eqref{eq:tran_bd}.
\end{remark}

\begin{remark}
The right hand side of Eqn~\eqref{eq:tran_bd2} can be modified as in standard IEF-PCM by using the identity
\[
	-\mathcal S_\kappa
\left(\partial_{\mathbf n}\psi_0\right)
	= (2\pi - \mathcal D_\kappa) \psi_0,
	\qquad\mbox{on }\Gamma,
\]
where $\mathcal D_\kappa$ is the corresponding double layer boundary operator, 
see \cite{cances1998new,cances1997new}.
This allows to obtain a right hand side that only depends on the potential and not on the field which subsequently leads to simpler expressions in the contribution to the Fock-matrix, if coupled to a quantum mechanical Hamiltonian within a polarizable embedding. 
\end{remark}
\black

\section{Strategy}\label{sect:strategy}
In this section, we introduce a global strategy for solving Eqs \eqref{eq:tran_prob}--\eqref{eq:tran_bd} that are derived from the LPB equation \eqref{eq:linpb}.
Then, we present how the domain decomposition method can be applied to solve the two PDEs defined on $\Omega$, taking advantage of its particular geometrical structure (i.e., a union of overlapping balls).
The scheme of this section is inspired by \cite[Section 4.2 and 5]{quan2018domain}, our previous work for the case of non-ionic solvent.

\subsection{Global strategy}
We propose the following iterative procedure for solving Eqs \eqref{eq:tran_prob}--\eqref{eq:tran_bd}: 
let $g^0$ defined on $\Gamma$ be an initial guess for the Dirichlet condition $\psi_{\rm e}|_{\Gamma}$ and set $k=1$.
\begin{itemize}
\item[\tt{[1]}] 
{Solve the following Dirichlet boundary problem for $\psi_{\rm r}^k$:
\begin{equation}
\begin{array}{r@{}l}
\left\{
\begin{aligned}
-\Delta\psi_{\rm r}^k&= 0 && \mbox{in } \Omega ,\\
\psi_{\rm r}^k & = g^{k-1} - \psi_{0} && \mbox{on }\Gamma,
\end{aligned}
\right. \label{eq:dd_laplace}
\end{array}
\end{equation}
and derive its Neumann boundary trace $\partial_{\mathbf n} \psi_{\rm r}^k$ on $\Gamma$.
}
\item[\tt{[2]}]
{Solve the following Dirichlet boundary problem for $\psi_{\rm e}^k$:
\begin{equation}
\begin{array}{r@{}l}
\left\{
\begin{aligned}
- \Delta\psi_{\rm e}^k + \kappa^2 \psi_{\rm e}^k& = 0 && \mbox{in } \Omega ,\\
\psi_{\rm e}^k & = g^{k-1} && \mbox{on }\Gamma,
\end{aligned}
\right.
\end{array}\label{eq:dd_screen}
\end{equation}
and derive similarly its Neumann boundary trace $\partial_{\mathbf n}\psi_{\rm e}^k$ on $\Gamma$.
}
\item[\tt{[3]}] Build the charge density $\sigma_{\rm e}^k = \displaystyle\partial_{\mathbf n}\psi_{\rm e}^k - \frac{\varepsilon_1}{\varepsilon_2}\partial_{\mathbf n}\left(\psi_0+\psi_{\rm r}^k \right)$ and compute a new Dirichlet condition $g^k = \mathcal {S}_{\kappa} \sigma_{\rm e}^k$.
\item[\tt{[4]}] Compute the contribution $E_k^{\rm s}$ to the solvation energy based on $\psi_{\rm r}^k$ at the $k$-th iteration, set $k \leftarrow k+1$, go back to Step {\tt{[1]}} and repeat until the increment of interaction $|E_k^{\rm s}-E_{k-1}^{\rm s}|$ becomes smaller than a given tolerance $\tt Tol \ll 1$.
\end{itemize}
\begin{remark}
In order to provide a suitable initial guess of $g^0$ (defined on $\Gamma$), we consider the (unrealistic) scenario where the whole space $\mathbb R^3$ is covered by the solvent medium. 
Then, the electrostatic potential $\psi$ in this case is given explicitly by 
\begin{equation}
\psi(\mathbf x) = \sum_{i = 1}^{M}\frac{4\pi q_i}{\varepsilon_2}\frac{\exp\left({-\kappa |\mathbf x-\mathbf x_i|}\right)}{|\mathbf x - \mathbf x_i|}, \quad \forall \mathbf x\in \mathbb R^3,
\end{equation}
see details in \cite[Section 1.3.2]{lamm2003poisson}. 
As a consequence, we choose $g^0$ as this potential restricted on $\Gamma$.
\end{remark}
\begin{remark}\label{rmk:no_ite}
The above global strategy is an iterative procedure, which is presented for an easier understanding. 
However, the final convergent solution satisfies, after discretization, a global linear system that can be solved by different linear solvers.
We will address this issue in the later Section \ref{sect:disc}.
\end{remark}

\subsection{Domain decomposition (DD) scheme}\label{subsect:dd}
The Schwarz's domain decomposition method \cite{quarteroni1999domain} is a good choice to solve the PDE defined on a complex domain which can be composed as a union of overlapping and possibly simple subdomains.
According to the definition of $\Omega$, we have a natural domain decomposition as follows
$$\Omega = \bigcup_{j=1}^{M}\Omega_j,\quad \Omega_j = B_{r_j}(\mathbf x_j),$$
where each $\Omega_j$ denotes the $j$-th VDW-ball with center $\mathbf x_j$ and radius $r_j$.
As a consequence, the Schwarz's domain decomposition method can be applied to solve the PDEs \eqref{eq:dd_laplace} and \eqref{eq:dd_screen}.

\begin{figure}
\centering
\includegraphics[width = 2 in]{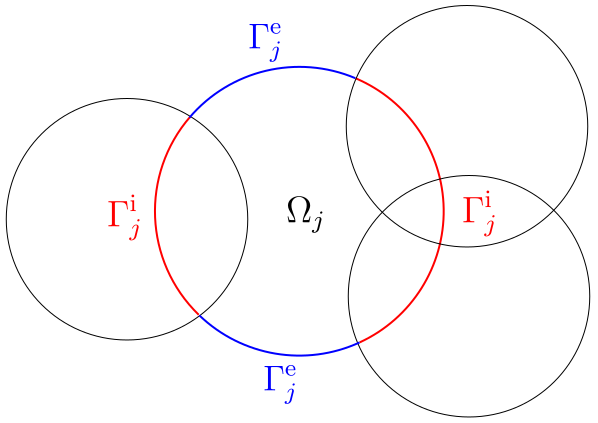}
\caption{2D schematic diagram of $\Gamma_j^{\rm i}$ (red) and $\Gamma_j^{\rm e}$ (blue) associated with $\Omega_j$.}\label{fig:dd}
\end{figure}

Similar to the ddCOSMO method \cite{cances2013domain}, Eq. \eqref{eq:dd_laplace} is equivalent to the following coupled local equations, each restricted to $\Omega_j$:
\begin{equation}
\begin{array}{r@{}l}
\left\{
\begin{aligned}
-\Delta\psi_{\rm r}|_{\Omega_j}&= 0 && \mbox{in } \Omega_j,\\
\psi_{\rm r}|_{\Gamma_j} & = \phi_{{\rm r},j} && \mbox{on }\Gamma_j,
\end{aligned}
\right.
\end{array}\label{eq:subeq}
\end{equation}
where $\Gamma_j = \partial \Omega_j$ and
\begin{equation}
\begin{array}{r@{}l}
\phi_{{\rm r},j} = 
\left\{
\begin{aligned}
& \psi_{\rm r} && \mbox{on }\Gamma_j^{\rm i},\\
&  g-\psi_0 && \mbox{on }\Gamma_{j}^{\rm e}.
\end{aligned}
\right.
\end{array}\label{eq:subeq_bdy}
\end{equation}
Here, we omit the superscript due to the (outer) iteration index $k$.
$\Gamma_j^{\rm e}$ is the external part of $\Gamma_j$ not contained in any other ball $\Omega_i$ ($i\neq j$), i.e., $\Gamma_j^{\rm e} = \Gamma\cap \Gamma_j$; $\Gamma_j^{\rm i}$ is the internal part of $\Gamma_j$, i.e., $\Gamma_j^{\rm i} = \Omega\cap \Gamma_j$ (see Figure \ref{fig:dd} for an illustration).
%
Similarly, Eq. \eqref{eq:dd_screen} is equivalent to the following coupled local equations, each restricted to $\Omega_j$:
\begin{equation}
\begin{array}{r@{}l}
\left\{
\begin{aligned}
-\Delta\psi_{\rm e}|_{\Omega_j} + \kappa^2 \psi_{\rm e}|_{\Omega_j}&= 0 && \mbox{in } \Omega_j,\\
\psi_{\rm e}|_{\Omega_j} & = \phi_{{\rm e},j} && \mbox{on }\Gamma_j,
\end{aligned}
\right.
\end{array}\label{eq:subeq2}
\end{equation}
where 
\begin{equation}
\begin{array}{r@{}l}
\phi_{{\rm e},j} = 
\left\{
\begin{aligned}
& \psi_{\rm e} && \mbox{on }\Gamma_j^{\rm i},\\
&  g&& \mbox{on }\Gamma_{j}^{\rm e}.
\end{aligned}
\right.
\end{array}\label{eq:subeq2_bdy}
\end{equation}

Note that the Dirichlet conditions that appear in \eqref{eq:subeq_bdy} and \eqref{eq:subeq2_bdy} are implicit since $\psi_{\rm r}$ (resp.  $\psi_{\rm e}$) is not known on $\Gamma_j^{\rm i}$. Hence, given the Dirichlet boundary condition on $\Gamma$, an iterative procedure must be applied to solve the coupled equations \eqref{eq:subeq}--\eqref{eq:subeq_bdy} (resp. \eqref{eq:subeq2}--\eqref{eq:subeq2_bdy}), such as the parallel Schwarz algorithm and the alternating Schwarz algorithm as presented in the ddCOSMO \cite{cances2013domain}.
For example, the idea of the parallel algorithm is to solve each local problem based on the boundary condition of the neighboring solutions derived from the previous iteration. 
In this iterative procedure, the computed value of $\psi_{\rm r}|_{\Gamma_j^{\rm i}}$ (resp. $\psi_{\rm e}|_{\Gamma_j^{\rm i}}$) is updated step by step and converges to the exact value.  

However, the parallel and the alternating Schwarz algorithms might not be the most efficient way to solve such a set of equations, but they are well-suited to illustrate the idea of the domain decomposition method. 
In fact, the global linear system obtained after discretization can be solved by different linear solvers (for example, the GMRES method). 
We will discuss more about this in Section \ref{sect:disc}.
Before that, we shall develop two single-domain solvers: a Laplace solver and a HSP solver in the unit ball.

\section{Single-domain solvers}\label{sect:solver}
In this section, we develop two single domain solvers in the unit ball respectively for solving Eq. \eqref{eq:subeq} and Eq. \eqref{eq:subeq2} within the domain decomposition scheme.

\subsection{Laplace solver}\label{sect:laplace_solver}
Developing a Laplace solver in a ball is not difficult and has already been presented in our previous work including the ddCOSMO \cite{cances2013domain}, the ddPCM \cite{stamm2016new} and the ddPCM-SES \cite{quan2018domain}. 
For the sake of completeness, we recall briefly the Laplace solver in the following content.

We want to solve the Laplace equation \eqref{eq:subeq} defined on $\Omega_j$.
Without loss of generality, we consider the following Laplace equation defined in the unit ball
\begin{equation}
\begin{array}{r@{}l}
\left\{
\begin{aligned}
-\Delta u_{\rm r} & = 0 && \mbox{in } B_1(\mathbf 0),\\
u_{\rm r} & = \phi_{\rm r} && \mbox{on } \mathbb S^2.
\end{aligned}
\right.
\end{array}\label{eq:laplace}
\end{equation}
Its unique solution in $H^1(B_1(\mathbf 0))$ can be written as
\begin{equation}\label{eq:sol_laplace}
u_{\rm r}(r,\theta,\varphi)  =\sum_{\ell=0}^\infty \sum_{m=-\ell}^\ell [\phi_{\rm r}]_\ell^m\, r^\ell \, Y_\ell^m(\theta,\varphi), \quad  0\leq r\leq 1,~0\leq \theta\leq \pi,~0\leq\varphi <2\pi.
\end{equation}
Here, $Y_\ell^m$ denotes the (real orthonormal) spherical harmonic of degree $\ell$ and order $m$ defined on $\mathbb S^2$
and $$ [\phi_{\rm r}]_\ell^m = \int_{\mathbb S^2} \phi_{\rm r}(\mathbf s) Y_\ell^m(\mathbf s) d\mathbf s,$$ is the real coefficient of $u_{\rm r}$ corresponding to the mode $Y_\ell^m$. 
Then, $u_{\rm r}$ can be numerically approximated by $\widetilde u_{\rm r}$ in the discretization space spanned by a truncated basis of spherical harmonics $\{Y_\ell^m\}_{0\leq \ell\leq \ell_{\rm max},~-\ell\leq m \leq \ell}$, defined as
\begin{equation}
\widetilde u_{\rm r}(r,\theta,\varphi)  =\sum_{\ell=0}^{\ell_{\rm max}}\sum_{m=-\ell}^\ell [\widetilde \phi_{\rm r}]_\ell^m\, r^\ell \, Y_\ell^m(\theta,\varphi), \quad 0\leq r\leq 1,~0\leq \theta\leq \pi,~0\leq\varphi <2\pi, \label{eq:num_laplace}
\end{equation}
where $\ell_{\rm max}$ denotes the maximum degree of spherical harmonics and 
\begin{equation}
[\widetilde \phi_{\rm r}]_\ell^m = \sum_{n = 1}^{N_{\rm leb}} w_n^{\rm leb} \phi_{\rm r}(\mathbf s_n) Y_\ell^m(\mathbf s_n).\label{eq:num_phi_r}
\end{equation}
Here, $\mathbf s_n\in \mathbb S^2$ represent Lebedev quadrature points \cite{0953-4075-40-23-004}, $w_n^{\rm leb}$ are the corresponding weights and $N_{\rm leb}$ is the number of Lebedev quadrature points.

\subsection{HSP solver}
\begin{figure}
\centering
\includegraphics[width = 5 in,clip,trim={1.5in 0 1.5in 0}]{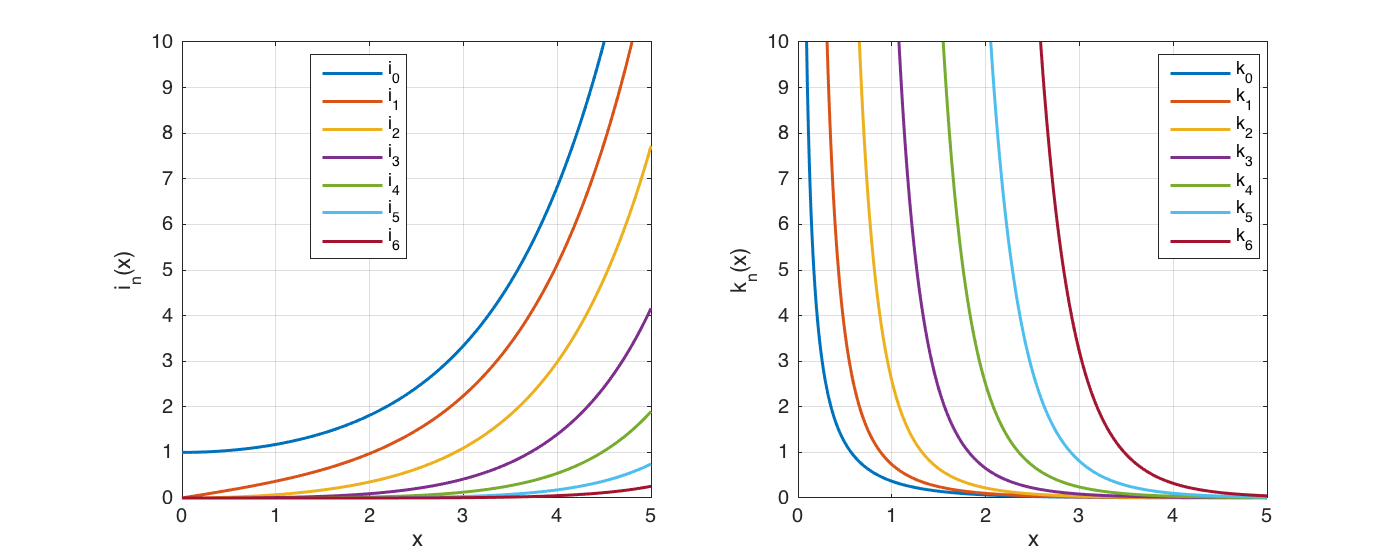}
\caption{The modified spherical Bessel functions of the first kind ($i_n$, left) and the second kind ($k_n$, right) with $\kappa  = 1$.}\label{fig:bessel_ms}
\end{figure}

We now want to solve the HSP equation \eqref{eq:subeq2} defined on $\Omega_j$.
Without loss of generality, we consider the following HSP equation defined in the unit ball
\begin{equation}
\begin{array}{r@{}l}
\left\{
\begin{aligned}
-\Delta u_{\rm e} + \kappa ^2  u^2_{\rm e}& = 0 && \mbox{in } B_1(\mathbf 0),\\
u_{\rm e} & = \phi_{\rm e} && \mbox{on } \mathbb S^2.
\end{aligned}
\right.
\end{array}\label{eq:screen}
\end{equation}
Solving the above HSP equation in spherical coordinates by separation of variables, the radial equation corresponding to the angular dependency $Y_\ell^m$ has the form
\begin{equation}
\frac{1}{R}\frac{d}{dr}\left(r^2\frac{dR}{dr}\right) = \kappa^2 r^2+\ell(\ell+1),\quad \ell \geq 0,
\end{equation}
that is,
\begin{equation}
r^2\frac{d^2R}{dr^2} + 2r \frac{dR}{dr}- (\kappa^2 r^2+\ell(\ell+1))R = 0,
\end{equation}
which is called the modified spherical Bessel equation \cite{arfken2012mathematical}.
This equation has two linearly independent solutions as follows
\begin{equation}
i_\ell( r) = \sqrt{\frac{\pi}{2\kappa r}} I_{\ell+\frac 1 2} (\kappa r),\quad k_\ell(r) = \sqrt{\frac{2}{\pi \kappa r}} K_{\ell+\frac 1 2} (\kappa r),
\end{equation}
where $i_\ell$ and $k_\ell$ are the modified spherical Bessel functions of the first and second kind associated with $\kappa$, see \cite[Chapter 14]{arfken2012mathematical} for details and Figure \ref{fig:bessel_ms} for an illustration.
Here, $I_{\alpha}(x)$ and $K_{\alpha}(x)$ with subscript $\alpha$ are the modified Bessel functions of the first and second kind \cite{abramowitz1964handbook}.
\begin{remark}
$I_{\alpha}(x)$ and $K_{\alpha}(x)$ satisfy the modified Bessel equation
\begin{equation}
x^2 \frac{d^2 f}{dx^2} + x \frac{df}{dx} - (x^2 + \alpha^2)f = 0.
\end{equation}
In fact, $I_\alpha$ and $K_\alpha$ are exponentially growing and decaying functions, respectively.
\end{remark}
Since $k_\ell\rightarrow \infty$ as $r\rightarrow 0$, we are interested in the family $i_\ell$ of the first kind.
That is, we write the solution to \eqref{eq:screen} in the form of 
\begin{equation}\label{eq:sol_screen}
u_{\rm e}(r,\theta,\varphi)  =\sum_{\ell=0}^\infty \sum_{m=-\ell}^\ell c_\ell^m\, i_\ell( r) \, Y_\ell^m(\theta,\varphi), \quad  0\leq r\leq 1,~0\leq \theta\leq \pi,~0\leq\varphi <2\pi,
\end{equation}
where $c_\ell^m $ is the coefficient of the mode $Y_{\ell}^m$.
With the same discretization as in Section \ref{sect:laplace_solver}, we derive the following approximate solution similar to Eq. \eqref{eq:num_laplace}:
\begin{equation}\label{eq:num_screen}
\widetilde u_{\rm e}(r,\theta,\varphi)  =\sum_{\ell=0}^{\ell_{\rm max}}\sum_{m=-\ell}^\ell [\widetilde \phi_{\rm e}]_\ell^m\,\frac{i_\ell( r) }{i_\ell (1)}\, Y_\ell^m(\theta,\varphi), \quad 0\leq r\leq 1,~0\leq \theta\leq \pi,~0\leq\varphi <2\pi,
\end{equation}
where $[\widetilde \phi_{\rm e}]_\ell^m$ is given similar to \eqref{eq:num_phi_r} as follows
\begin{equation}
[\widetilde \phi_{\rm e}]_\ell^m = \sum_{n = 1}^{N_{\rm leb}} w_n^{\rm leb} \phi_{\rm e}(\mathbf s_n) Y_\ell^m(\mathbf s_n),\label{eq:num_phi_e}
\end{equation}
with the same notations $w_n^{\rm leb}$ and $N_{\rm leb}$ as above.

\begin{remark}
According to the fact that (see \cite[page 708]{arfken2012mathematical})
$$
i_\ell(r) \approx \frac{(\kappa r)^\ell } {(2\ell+1)!!},\quad \mbox{when $r> 0$ is very small,}
$$
we have $\frac{i_\ell(r)}{i_\ell(1)} \rightarrow r^\ell$ as $\kappa \rightarrow 0$.
Therefore, if $\phi_{\rm r} = \phi_{\rm e}$,  then $\widetilde u_{\rm e}\rightarrow \widetilde u_{\rm r}$ as $\kappa \rightarrow 0$.
This means that the solution to Eq. \eqref{eq:screen} tends to the solution to Eq. \eqref{eq:laplace}, which makes sense.
\end{remark}


\section{Discretization}\label{sect:reform_disc}
The global strategy in Section \ref{sect:strategy} in combination with the domain decomposition schemes for solving Eq. \eqref{eq:dd_laplace} and Eq. \eqref{eq:dd_screen} is an iterative procedure.
This implies that the proposed algorithm can be parallelized since only a group of local problems on each $\Omega_j$ are solved.
However, as mentioned in Remark \ref{rmk:no_ite}, we will solve the problem in a global way, meaning that we finally solve a global linear system derived from discretization.
To present this, we first introduce a reformulation of the coupling conditions and then present the global linear system for its discretization.

\subsection{Reformulation}
Let $\chi_i$ be the characteristic function of $\Omega_i$, i.e.,
\begin{equation}\label{eq:chi}
\begin{array}{r@{}l}
\chi_i(\mathbf x) \coloneqq \left\{
\begin{aligned}
&1 &&\mbox{if }\mathbf x\in \Omega_i \\
&0 &&\mbox{if }\mathbf x\not\in \Omega_i
\end{aligned}
\right.
\end{array}
\end{equation}
and then let
\begin{equation}\label{eq:w_ji}
w_{ji}(\mathbf x) \coloneqq \frac{\chi_i(\mathbf x)}{|\mathcal N(j,\mathbf x)|} = \frac{\chi_i(\mathbf x)}{\sum_{i\neq j}\chi_i(\mathbf x)},\quad \forall \mathbf x\in \Gamma_j,
\end{equation}
where $\mathcal N(j,\mathbf x)$ represents the index set of all balls containing $\mathbf x$.
Here, we make the convention that in the case of $\left|\mathcal N(j,\mathbf x)\right| = 0$ (i.e., $\mathbf x\in \Gamma_j^{\rm e}$), we define $w_{ji}(\mathbf x) = 0,~\forall i$.
Furthermore, $\forall\mathbf x\in \Gamma_j$, we define
\begin{equation}
\begin{array}{r@{}l}
\chi_j^{\rm e}(\mathbf x) \coloneqq \left\{
\begin{aligned}
&1 &&\mbox{if }\mathbf x\in \Gamma_j^{\rm e} ,\\
&0 &&\mbox{if } \mathbf x\in \Gamma_j^{\rm i} ,
\end{aligned}
\right.
\end{array}
\end{equation}
which is equivalent to
\begin{equation}
\chi_j^{\rm e}(\mathbf x) = 1 - \sum_{i\neq j} w_{ji}(\mathbf x), \quad \forall \mathbf x\in \Gamma_j.
\end{equation}

\subsubsection{DD scheme}
There are two local coupling conditions in the DD scheme, i.e., Eq. \eqref{eq:subeq_bdy} and \eqref{eq:subeq2_bdy}, respectively for coupling the local Laplace equation \eqref{eq:subeq} and the local HSP equation \eqref{eq:subeq2}.
Based on the above-defined notations, Eq. \eqref{eq:subeq_bdy} can be recast as
\begin{equation}\label{eq:psi_r_j}
\psi_{\rm r}|_{\Gamma_j}(\mathbf x) -\sum_{i\neq j} w_{ji}(\mathbf x)\,\psi_{\rm r}|_{\Omega_i}(\mathbf x)=\chi_j^{\rm e}(\mathbf x) \left( g(\mathbf x)-\psi_0(\mathbf x) \right),\quad \forall \mathbf x\in\Gamma_j.
\end{equation}
Similarly, Eq. \eqref{eq:subeq2_bdy} can be recast as
\begin{equation}\label{eq:psi_e_j}
\psi_{\rm e}|_{\Gamma_j}(\mathbf x) -\sum_{i\neq j} w_{ji}(\mathbf x)\,\psi_{\rm e}|_{\Omega_i}(\mathbf x)=\chi_j^{\rm e}(\mathbf x) g(\mathbf x),\quad \forall \mathbf x\in\Gamma_j.
\end{equation}

\subsubsection{Boundary coupling condition}
There is a global coupling condition on $\Gamma$, i.e., Eq. \eqref{eq:tran_bd2}, between the global Laplace equation and the global HSP equation defined in $\Omega$ (see \eqref{eq:tran_prob}), involving the nonlocal operator $\mathcal S_{\kappa}$:
\begin{equation}\label{eq:g}
g(\mathbf x) = \psi_{\rm e}|_{\Gamma}(\mathbf x) = \mathcal S_\kappa\left(\partial_{\mathbf n}\psi_{\rm e} - \frac{\varepsilon_1}{\varepsilon_2}\partial_{\mathbf n}\left(\psi_0+\psi_{\rm r} \right)\right) (\mathbf x),\quad \forall\mathbf x\in \Gamma.
\end{equation}
The single-layer operator $\mathcal S_\kappa$ involves an integral over the whole solute-solvent boundary $\Gamma$ which seems difficult to compute at a first glance.
\black
We introduce a technique to compute the integral of $\mathcal S_\kappa$ efficiently.
For each sphere $\Gamma_i$, we define a local single-layer potential $\widetilde{\mathcal S}_{\kappa,\Gamma_i}$ as follows
\begin{equation}
\widetilde{\mathcal S}_{\kappa,\Gamma_i}\sigma (\mathbf x)  \coloneqq  \int_{\Gamma_i} \frac{\exp\left({-\kappa |\mathbf x-\mathbf y|}\right)\sigma(\mathbf y)}{4\pi |\mathbf x - \mathbf y|}, \quad \forall \mathbf x\in \mathbb R^3,
\end{equation}
where $\sigma$ is an arbitrary function in $H^{-\frac 1 2}(\Gamma_i)$.
As a consequence, $\forall \sigma\in H^{-\frac 1 2}(\Gamma)$, we have
\begin{equation}
\mathcal S_\kappa \sigma  = \sum_{i = 1}^M \widetilde{\mathcal S}_{\kappa,\Gamma_i}\left(\chi_i^{\rm e} \sigma\right),
\end{equation}
where $\chi_i^{\rm e} \sigma$ extends $\sigma|_{\Gamma_i^{\rm e}}$ by zero to the whole sphere $\Gamma_i$.
The above equation implies that the integral over $\Gamma$ can be divided into a group of integrals respectively over each sphere $\Gamma_i$.
Therefore, Eq. \eqref{eq:g} can be recast as
\begin{equation}\label{eq:g_tran_form}
\begin{array}{r@{}l}
\begin{aligned}
g(\mathbf x) &= \sum_{i = 1}^{M} \widetilde{\mathcal S}_{\kappa,\Gamma_i} \left[\chi_i^{\rm e} \left(\partial_{\mathbf n}\psi_{\rm e} - \frac{\varepsilon_1}{\varepsilon_2}\partial_{\mathbf n}\left(\psi_0+\psi_{\rm r} \right)\right) \right] (\mathbf x), \quad \forall \mathbf x\in \Gamma,
\end{aligned}
\end{array}
\end{equation}
{which is used for updating the boundary potential in the global strategy in Section \ref{sect:strategy}.}

\subsection{Linear system}\label{sect:disc} 
In this part, we first present the discretization of the above reformulation and then introduce the global linear system derived from this discretization.

\subsubsection{Local truncation in balls}
In Section \ref{sect:solver}, without the loss of generalization, we have presented the discretization of the solutions to the Laplace equation and the HSP equation defined in the unit ball.

Based on this, for each sphere $\Gamma_j$, we first approximate $\psi_{\rm r}|_{\Gamma_j}$ and $\psi_{\rm e}|_{\Gamma_j}$ respectively by a linear combination of spherical harmonics $\{Y_\ell^m\}$ with $0\leq\ell\leq\ell_{\rm max}$ and $-\ell\leq m\leq\ell$ as follows
\begin{equation}\label{eq:psi_r_j_num}
\psi_{\rm r} |_{\Gamma_j}( \mathbf x_j + r_j\mathbf s) = \sum_{\ell=0}^{\ell_{\rm max}}\sum_{m=-\ell}^\ell [X_{\rm r}]_{j \ell m}\, Y_\ell^m(\mathbf s), \quad \mathbf s\in\mathbb S^2,
\end{equation}
and
\begin{equation}
\psi_{\rm e} |_{\Gamma_j}( \mathbf x_j + r_j\mathbf s) = \sum_{\ell=0}^{\ell_{\rm max}}\sum_{m=-\ell}^\ell [X_{\rm e}]_{j\ell m} Y_\ell^m(\mathbf s), \quad \mathbf s\in\mathbb S^2,
\end{equation}
where $[X_{\rm r}]_{j\ell m}$ and $[X_{\rm e}]_{j\ell m}$ are unknown coefficients of the mode $Y_\ell^m$ respectively associated with $\psi_{\rm r}|_{\Gamma_j}$ and $\psi_{\rm e}|_{\Gamma_j}$.
Here, for any point $\mathbf x\in \Gamma_j$, we actually use its spherical coordinates $(r_j,\mathbf s)$ s.t. $\mathbf x = \mathbf x_j + r_j\mathbf s$.
According to the Laplace solver and the HSP solver presented in Section \ref{sect:solver}, we deduce directly
\begin{equation}\label{eq:psi_r_omega_i}
\psi_{\rm r} |_{\Omega_i}(\mathbf x_i + r\mathbf s) = \sum_{\ell'=0}^{\ell_{\rm max}}\sum_{m'=-\ell'}^{\ell'} [X_{\rm r}]_{i\ell' m'} \,\left(\frac{r}{r_i}\right)^{\ell'}\, Y_{\ell'}^{m'}(\mathbf s), \quad 0\leq r\leq r_i,~\mathbf s\in \mathbb S^2,
\end{equation}
and
\begin{equation}
\psi_{\rm e} |_{\Omega_i}(\mathbf x_i + r\mathbf s) = \sum_{\ell'=0}^{\ell_{\rm max}}\sum_{m'=-\ell'}^{\ell'} [X_{\rm e}]_{i\ell' m'} \,\frac{i_{\ell'}\left({r}\right)}{i_{\ell'}\left({r_i}\right)} \,Y_{\ell'}^{m'}(\mathbf s), \quad 0\leq r\leq r_i,~\mathbf s\in \mathbb S^2,
\end{equation}
where for any point $\mathbf x\in \Omega_i$, we take its spherical coordinates $(r,\mathbf s)$ s.t. $\mathbf x = \mathbf x_i + r\mathbf s$.
Also, for each sphere $\Gamma_i$, we can compute the normal derivative of $\psi_{\rm r}$ on $\Gamma_i^{\rm e}$ as follows
\begin{equation}\label{eq:partial_psi_r}
\partial_{\mathbf n}\psi_{\rm r}(\mathbf x_i + r_i\mathbf s) =  \sum_{\ell'=0}^{\ell_{\rm max}}\sum_{m'=-\ell'}^{\ell'} [X_{\rm r}]_{i\ell' m'} \,\left(\frac{\ell'}{r_i}\right)\, Y_{\ell'}^{m'}(\mathbf s),\quad \mathbf x_i + r_i \mathbf s\in \Gamma_i^{\rm e},
\end{equation}
and the normal derivative of $\psi_{\rm e}$ on $\Gamma_i^{\rm e}$ 
\begin{equation}\label{eq:partial_psi_e}
\partial_{\mathbf n}\psi_{\rm e}(\mathbf x_i + r_i\mathbf s) =  \sum_{\ell'=0}^{\ell_{\rm max}}\sum_{m'=-\ell'}^{\ell'} [X_{\rm e}]_{i \ell' m'} \,\frac{i_{\ell'}'\left({ r_i}\right)}{i_{\ell'}\left({}{r_i}\right)}\, Y_{\ell'}^{m'}(\mathbf s),\quad \mathbf x_i + r_i \mathbf s \in \Gamma_i^{\rm e},
\end{equation}
where $i'_{\ell'}$ represents the derivative of $i_{\ell'}$.
\begin{remark}
We compute $i'_{\ell'}(r_i)$ according to the following equation (see \cite[page 707]{arfken2012mathematical} for the derivation)
\begin{equation}
(2n+1) i'_n(r_i) = n r_i i_{n-1}(r_i)+(n+1) r_i i_{n+1}(r_i),
\end{equation}
and in analogy, we compute $k'_\ell(r_i)$ used in the Appendix \ref{append:1} as follows
\begin{equation}
-(2n+1) k'_n(r_i) = n r_i k_{n-1}(r_i)+(n+1) r_i k_{n+1}(r_i).
\end{equation}
\end{remark}

\subsubsection{Discretization}\label{sect:disc_couplingconds}

So far, we have written the Ansatz for the unknowns $\psi_{\rm r}|_{\Gamma_j}$, $\psi_{\rm e}|_{\Gamma_j}$, respectively $\psi_{\rm r}|_{\Omega_i}$, $\psi_{\rm e}|_{\Omega_i}$ with normal derivatives $\partial_{\mathbf n}\psi_{\rm r}$ and $\partial_{\mathbf n}\psi_{\rm e}$, following Eqs \eqref{eq:psi_r_j_num} -- \eqref{eq:partial_psi_e}, which depend on the unknowns $X_{\rm r}$ and $X_{\rm e}$.
This allows us to discretize the coupling conditions \eqref{eq:psi_r_j}, \eqref{eq:psi_e_j} and \eqref{eq:g_tran_form}, to derive a final linear system.

First, we replace the variable $\mathbf x\in \Gamma_j$ of Eq. \eqref{eq:psi_r_j} by $\mathbf x = \mathbf x_j+r_j\mathbf s$ with $\mathbf s\in \mathbb S^2$ and derive the equation for each sphere $\Gamma_j$ as follows
\begin{equation}
\begin{array}{r@{}l}
\begin{aligned}
&\psi_{\rm r}|_{\Gamma_j}(\mathbf x_j + r_j\mathbf s) -\sum_{i\neq j} w_{ji}(\mathbf x_j +r_j\mathbf s)\,\psi_{\rm r}|_{\Omega_i}(\mathbf x_j +r_j\mathbf s) \\
= & \chi_j^{\rm e}(\mathbf x_j +r_j\mathbf s)\left( g(\mathbf x_j +r_j\mathbf s)- \psi_0(\mathbf x_j +r_j\mathbf s)\right),
\end{aligned}
\end{array}
\end{equation}
which induces the following local equation by multiplying by $Y_\ell^m$ and integrating over $\mathbb S^2$ on both sides, $\forall j,\ell,m,$
\begin{equation}\label{eq:vf}
\begin{array}{r@{}l}
\begin{aligned}
&\left \langle \psi_{\rm r}|_{\Gamma_j}(\mathbf x_j + r_j\boldsymbol{\cdot}) -\sum_{i\neq j} w_{ji}(\mathbf x_j +r_j\boldsymbol{\cdot})\,\psi_{\rm r}|_{\Omega_i}(\mathbf x_j +r_j\boldsymbol{\cdot}), Y_\ell^m(\boldsymbol{\cdot})\right\rangle_{\mathbb S^2}\\ 
& = \left\langle \chi_j^{\rm e}(\mathbf x_j +r_j\boldsymbol{\cdot})\left( g(\mathbf x_j +r_j\boldsymbol{\cdot})- \psi_0(\mathbf x_j +r_j\boldsymbol{\cdot})\right),Y_{\ell}^m(\boldsymbol{\cdot}) \right\rangle_{\mathbb S^2}.
\end{aligned}
\end{array}
\end{equation}
Here, $\langle\cdot,\cdot\rangle_{\mathbb S^2}$ represents the integral over the unit sphere $\mathbb S^2$, which is numerically approximated using the Lebedev quadrature rule with $N_{\rm leb}$ points.
We therefore denote such a numerical integration over $\mathbb S^2$ by the notation $\langle \cdot,\cdot \rangle_{\mathbb S^2,N_{\rm leb}}$.
Eq. \eqref{eq:vf} can be rewritten in the form of a linear system
\begin{equation}\label{eq:linsys_A}
[\mathbf A X_{\rm r}]_{j \ell m}  = [G_X]_{j \ell m} + [G_0]_{j \ell m},\quad \forall j,\ell,m.
\end{equation}
Here,  $\mathbf A$ is a square matrix of dimension $M(\ell_{\rm max}+1)^2\times M(\ell_{\rm max}+1)^2$ and the $j\ell m$-th row of $\mathbf A X_{\rm r}$ is given by substituting \eqref{eq:psi_r_j_num} and \eqref{eq:psi_r_omega_i} into \eqref{eq:vf} as follows
\begin{equation}
\begin{array}{r@{}l}
\begin{aligned}
[\mathbf A X_{\rm r}]_{j \ell m} &  = \left\langle \psi_{\rm r}|_{\Gamma_j}(\mathbf x_j + r_j\boldsymbol{\cdot}) -\sum_{i\neq j} w_{ji}(\mathbf x_j +r_j\boldsymbol{\cdot})\,\psi_{\rm r}|_{\Omega_i}(\mathbf x_j +r_j\boldsymbol{\cdot}), Y_\ell^m(\boldsymbol{\cdot}) \right\rangle_{\mathbb S^2,N_{\rm leb}} \\
& =  [X_{\rm r}]_{j \ell m} - \sum_{i\neq j} \sum_{\ell', m'}\\
&\quad \left(\sum_{n = 1}^{N_{\rm leb}}  w_{n}^{\rm leb}w_{ji}(\mathbf x_j +r_j\mathbf s_n)\, \left(\frac{r_{ijn}}{r_i}\right)^{\ell'} Y_{\ell'}^{m'}(\mathbf s_{ijn})Y_\ell^m(\mathbf s_n)\right) [X_{\rm r}]_{i \ell' m'} ,
\end{aligned}
\end{array}
\end{equation}
where $(r_{ijn},\mathbf s_{ijn})$ is the spherical coordinate associated with $\Gamma_i$ of the point $\mathbf x_j +r_j \mathbf s_n$ s.t. 
$$ \mathbf x_j +r_j \mathbf s_n = \mathbf x_i + r_{ijn}\mathbf s_{ijn},\quad\mbox{with}\quad \mathbf s_{ijn} \in \mathbb S^2.$$
Furthermore, the $j\ell m$-th element of the column vector $G_X$ is given as
\begin{equation}
\begin{array}{r@{}l}
\begin{aligned}
[G_X]_{j\ell m} &= \left\langle \chi_j^{\rm e}(\mathbf x_j +r_j\boldsymbol{\cdot}) g(\mathbf x_j +r_j\boldsymbol{\cdot}) ,Y_{\ell}^m(\boldsymbol{\cdot})\right\rangle_{\mathbb S^2, N_{\rm leb}} \\
&=\sum_{n = 1}^{N_{\rm leb}} w_n^{\rm leb} \chi_j^{\rm e}(\mathbf x_j +r_j\mathbf s_n) g(\mathbf x_j +r_j\mathbf s_n) Y_{\ell}^m(\mathbf s_n),
\end{aligned}
\end{array}\label{eq:G_jlm}
\end{equation}
which depends on the unknowns $X_{\rm r}$ and $X_{\rm e}$ through $g$ given by Eq. \eqref{eq:g_tran_form}.
The notation $X$ denotes the column of all unknowns, i.e.,
\begin{equation}
X = \left(\begin{array}{c} X_{\rm r} \\ X_{\rm e} \end{array}\right) \in\mathbb R^{2M(\ell_{\rm max}+1)^2}.
\end{equation}
Similarly, the $j\ell m$-th element of the column vector $G_0$ is given as
\begin{equation}
\begin{array}{r@{}l}
\begin{aligned}
[G_0]_{j\ell m} = - \sum_{n = 1}^{N_{\rm leb}} w_n^{\rm leb} \chi_j^{\rm e}(\mathbf x_j +r_j\mathbf s_n) \psi_0(\mathbf x_j +r_j\mathbf s_n) Y_{\ell}^m(\mathbf s_n),
\end{aligned}
\end{array}
\end{equation}
which can be computed a priori, since it is independent of $X$.
\black

Similar to the linear system \eqref{eq:linsys_A} and according to Eq. \eqref{eq:psi_e_j} for each sphere $\Gamma_j$, we have another linear system in the form of matrices
\begin{equation}\label{eq:linsys_B}
[\mathbf B X_{\rm e}]_{j \ell m}  = [G_X]_{j \ell m}, \quad \forall j,\ell,m,
\end{equation}
where the square matrix $\mathbf B$ satisfies
\begin{equation}
\begin{array}{r@{}l}
\begin{aligned}
[\mathbf B  X_{\rm e}]_{j \ell m} &= [X_{\rm e}]_{j \ell m}- \sum_{i\neq j} \sum_{\ell', m'} \\
&\quad \left(\sum_{n = 1}^{N_{\rm leb}}  w_{n}^{\rm leb}\, w_{ji}(\mathbf x_j +r_j\mathbf s_n)\, \frac{i_{\ell'}\left({ r_{ijn}}\right)}{i_{\ell'}\left({}{r_i}\right)}\, Y_{\ell'}^{m'}(\mathbf s_{ijn})\, Y_\ell^m(\mathbf s_n) \right) [X_{\rm e}]_{i \ell' m'},
\end{aligned}
\end{array}
\end{equation}
and $[G_X]_{j \ell m}$ is given by \eqref{eq:G_jlm}.

So far, we have derived two linear systems of the form 
\begin{equation}\label{eq:linAB}
\begin{array}{r@{}l}
\left\{
\begin{aligned}
&\mathbf A \, X_{\rm r}  = G_X + G_0,\\
&\mathbf B \, X_{\rm e} = G_X,
\end{aligned}
\right.
\end{array}
\end{equation}
where $X_{\rm r}$ and $X_{\rm e}$ are the column vectors of unknowns $[X_{\rm r}]_{j\ell m}$ and $[X_{\rm e}]_{j\ell m}$ (respectively associated with the potentials $\psi_{\rm r}$ and $\psi_{\rm e}$). 
However, the column vector $G_X$ depending on both $X_{\rm r}$ and $X_{\rm e}$ is not specified yet.
To do this, the coupling condition \eqref{eq:g_tran_form} in terms of $g$ should be used (which has not been used yet).
Combining Eq. \eqref{eq:g_tran_form}  with \eqref{eq:G_jlm}, we deduce the following form of $G_X$,
\begin{equation}\label{eq:G}
G_X = F_0 -\mathbf C_1 X_{\rm r} -\mathbf C_2 X_{\rm e},
\end{equation}
where the column vector $F_0$ is associated with $\partial_{\mathbf n}\psi_0$, two dense square matrices $\mathbf C_1$ and $\mathbf C_2$ are respectively associated with $\partial_{\mathbf n}\psi_{\rm r}$ and $\partial_{\mathbf n}\psi_{\rm e}$.
Considering the complexity of the formulas of $F_0,~\mathbf C_1,~\mathbf C_2$, we present them in Append. \ref{append:1}.

\begin{remark}
The number of the intersection of one atom with others is bounded from above.
From the definition \eqref{eq:w_ji} of $w_{ji}$, we have $w_{ji}(\mathbf x_j +r_j\mathbf s_n) = 0$ if $r_{ijn} \geq r_i$.
Therefore, $\mathbf A$ and $\mathbf B$ are both sparse matrices for a large molecule.
\end{remark}

\subsection{Linear solver}\label{subsect:linsolv}
We finally obtain a global linear system written as 
\begin{equation}\label{eq:linsys0}
\left(\begin{array}{cc}\mathbf A+\mathbf C_1 & \mathbf C_2 \\ \mathbf C_1 & \mathbf B+\mathbf C_2\end{array}\right) \left(\begin{array}{c} X_{\rm r} \\ X_{\rm e} \end{array}\right)= \left(\begin{array}{c} G_0 + F_0 \\ F_0 \end{array}\right),
\end{equation}
where both $\mathbf A$ and $\mathbf B$ are sparse for a large molecule, but not $\mathbf C_1$ nor $\mathbf C_2$.

To solve this linear system \eqref{eq:linsys0}, the LU factorization method and the GMRES method can be used directly \cite{saad2003iterative,golub2012matrix}, where the first one gives an exact solution and the second one gives an approximate solution.
However, the global strategy introduced in Section \ref{sect:strategy} provides another iterative strategy as follows
\begin{equation}\label{eq:lin_scheme}
\left(\begin{array}{cc}\mathbf A & \mathbf 0 \\ \mathbf 0 & \mathbf B\end{array}\right) \left(\begin{array}{c} X_{\rm r}^k \\ X_{\rm e}^k \end{array}\right)= -\left(\begin{array}{cc}\mathbf C_1 & \mathbf C_2 \\ \mathbf C_1 & \mathbf C_2\end{array}\right) \left(\begin{array}{c} X_{\rm r}^{k-1} \\ X_{\rm e}^{k-1} \end{array}\right) + \left(\begin{array}{c} G_0 + F_0 \\ F_0 \end{array}\right),
\end{equation}
where $k$ denotes the (outer) iteration number as in Section \ref{sect:strategy},  $X_{\rm r}^k$ and $X_{\rm e}^k$ are respectively the values of $X_{\rm r}$ and $X_{\rm e}$ computed at the $k$-th iteration.
At the $k$-th iteration, we first update the right-hand side of Eq. \eqref{eq:lin_scheme}, based on the previously-computed  $X_{\rm r}^{k-1}$ and $X_{\rm e}^{k-1}$.
Then, we use the GMRES method to solve the linear system associated with $X_{\rm r}^{k}$ and $X_{\rm e}^{k}$.
To distinguish the GMRES iterations, we call the above iteration with index $k$ as the outer iteration.

\black

\section{Numerical results}\label{sect:numerical}
For an implicit solvation model, one important issue is to compute the solute-solvent interaction energy, to which the electrostatic contribution plays an important role.
In fact, the electrostatic solvation energy $E^{\rm s}$ is computed from the reaction potential $\psi_{\rm r}$ according to the following formula (see \cite{fogolari2002poisson} for the derivation of this formula)
\begin{equation}\label{eq:Es}
E^{\rm s} =  \frac 1 2 \int_{\mathbb R^3}\rho_{\rm M}(\mathbf r)\, \psi_{\rm r}(\mathbf r) \, d\mathbf r,
\end{equation}
where the solute's charge density $\rho_{\rm M}$ is given in Eq. \eqref{eq:rho} and $\psi_{\rm r}$ is obtained by solving the LPB equation.
The outer iteration stops if ${\tt inc}_k<\tt Tol$, where ${\tt Tol}$ is the stopping tolerance and  
\begin{equation}
{\tt inc}_k \coloneqq \frac{|E^{\rm s}_{k}-E^{\rm s}_{k-1}|}{|E^{\rm s}_k|},
\end{equation}
where $E^{\rm s}_{k}$ denotes the electrostatic solvation energy computed at the $k$-th iteration.
In the following content, we study the electrostatic solvation energy computed numerically by the ddLPB, which has been implemented in both Matlab and Fortran.
Our Fortran code is based on the ddCOSMO and ddPCM codes written by Lipparini et al. (see GitHub link \cite{lipparini2015}).

By default, we take the dielectric permittivity in the solute cavity as in vacuum, that is, $\varepsilon_1 = 1$, and take the solvent to be water with the dielectric permittivity $\varepsilon_2 = 78.54$ at the room temperature $T = 298.15$K ($25^\circ$C).
Further, we set the Debye-H\"uckel screening constant to $\kappa=0.1040~ {\text \AA}^{-1}$, for an ionic strength $I = 0.1$ {\rm molar}.
The solute cavity is chosen as the VDW-cavity.
The atomic centers, charges and VDW radii are obtained from the PDB files \cite{berman2000protein} and the PDB2PQR package \cite{dolinsky2004pdb2pqr,dolinsky2007pdb2pqr} with the PEOEPB force field.
By default, the stopping tolerance introduced in Section \ref{subsect:linsolv} is set to ${\tt Tol} = 10^{-4}$, while the GMRES tolerance in the ddLPB is set to $10^{-8}$.
In the following content, these default parameters are used if they are not specified.

\subsection{Kirkwood model}
To test the ddLPB, we start from the Kirkwood model which has the explicit analytical solution, see  \cite{kirkwood1934theory,geng2015boundary,nguyen2017accurate}. 
In this model, there is only one sphere but with multiple charges distributed in this sphere.
We consider the following six cases with the sphere radii all set to $2$\AA.
\begin{itemize}
\item	Case 0 (Born model). One positive unit charge placed at $(0,0,0)$.
\item Case 1. Two positive unit charges placed at $(1, 0, 0)$ and $(-1, 0, 0)$.
\item Case 2. Two positive unit charges placed at $(1, 0, 0)$ and $(-1, 0, 0)$, and two negative unit charges placed at $(0, 1, 0)$ and $(0,-1, 0)$.
\item Case 3. Two positive unit charges placed at $(1.2, 0, 0)$ and $(-1.2, 0, 0)$, and two negative unit charges
symmetrically placed at $(0, 1.2, 0)$ and $(0,-1.2, 0)$.
\item Case 4. Six Positive unit charges placed at $(0.4, 0, 0)$, $(0, 0.8, 0)$,  $(0, 0, 1.2)$, $(0, 0,-0.4)$, $(-0.8, 0, 0)$ and $(0,-1.2, 0)$.
\item Case 5. Six positive unit charges placed at $(0.2,0.2, 0.2)$, $(0.5, 0.5, 0.5)$, \linebreak$(0.8, 0.8,0.8)$, $(-0.2, 0.2,-0.2)$, $(0.5,-0.5, 0.5)$ and $(-0.8,-0.8,-0.8)$.
\end{itemize}
The first case is also called the Born model and the other fives cases are recommended in \cite{nguyen2017accurate} to be tested for PB solvers.
Table \ref{tab:kirkwood} lists the electrostatic solvation energies computed by the ddLPB as well as the relative errors, which are computed according to the following definition 
\begin{equation}
\mbox{\tt RE} \coloneqq \frac{|E^{\rm s} - E^{\rm s}_{\rm ex}|}{|E^{\rm s}_{\rm ex}|}.
\end{equation}
Here, $E^{\rm s}$ denotes the electrostatic solvation energy computed by the ddLPB, while $E^{\rm s}_{\rm ex}$ denotes the exact result.
In the table, for a fixed $\ell_{\rm max}$, the corresponding $N_{\rm leb}$ ensures the accuracy of computing the scalar products of two arbitrary spherical harmonics in the function basis.

\begin{table}
\caption{Electrostatic solvation energies ({\rm kcal}/{\rm mol}) of different Kirkwood models computed by the ddLPB with different discretization parameters $\ell_{\rm max}$ and $N_{\rm leb}$, where $\kappa = 0 $. {\tt RE} represents the relative error of the ddLPB result with respective to the exact result.}\label{tab:kirkwood}
\scriptsize
\centering
\scalebox{1.}{
\renewcommand\arraystretch{1.2}
\begin{tabular}{c ccccccccc}
\toprule
 & & \multicolumn{2}{c}{Case 0} &  \multicolumn{2}{c}{Case 1} &  \multicolumn{2}{c}{Case 2}   \\ \cline{3-8}
{\bfseries $\ell_{\rm max}$}& {\bfseries $N_{\rm leb}$} & $E^{\rm s}$ &{\tt RE}&$E^{\rm s}$ & {\tt RE}& $E^{\rm s}$ &{\tt RE} \\
\midrule
3 & 26 & -81.9589 & 0 &  -349.5532 & 1.3762e-04 & -64.0454  &  2.0606e-02 \\
\hline
5 & 50 & -81.9589 & 0 &  -349.4132  & 2.6294e-04 &-62.5341   &  3.4772e-03 \\
\hline
7 &86 & -81.9589 & 0 &  -349.5064  & 3.7195e-06 & -62.7587&   1.0199e-04\\
\hline
9 & 146 & -81.9589 & 0 &  -349.5043  & 2.2890e-06 &-62.7508  &  2.3904e-05\\
\hline
11 & 194&  -81.9589 & 0 &  -349.5052 & 2.8612e-07 &-62.7524  & 1.5936e-06 \\
\hline
 \multicolumn{2}{c}{Exact} &  -81.9589 & &  -349.5051  & &  -62.7523   &\\
\bottomrule
\toprule
 & & \multicolumn{2}{c}{Case 3}  & \multicolumn{2}{c}{Case 4} & \multicolumn{2}{c}{Case 5}  \\\cline{3-8} 
{\bfseries $\ell_{\rm max}$} & {\bfseries $N_{\rm leb}$} & $E^{\rm s}$ & {\tt RE}& $E^{\rm s}$&{\tt RE} & $E^{\rm s}$& {\tt RE}\\
\midrule  
3 & 26 & -141.3629  &  4.5417e-02 &-2991.4727  &  9.8768e-04 &-3114.9078   & 2.7422e-03  \\
\hline
5 & 50  & -133.2890  &  1.4292e-02 &-2988.0565  & 1.5543e-04 &-3126.2105  & 8.7643e-04  \\
\hline
7& 86 &  -135.3437  &  9.0296e-04 &-2988.5869  &   2.2051e-05 &-3124.0588  & 1.8755e-04 \\
\hline
9 & 146 &  -135.1534  &  5.0436e-04 &-2988.4893  &   1.0607e-05 &-3123.5037   & 9.8288e-06\\
\hline
11 & 194 &  -135.2189 &  1.9967e-05 &-2988.5196  &   4.6846e-07 &-3123.5193   & 1.4823e-05 \\
\hline
 \multicolumn{2}{c}{Exact} &  -135.2216  & &  -2988.5210  &  &  -3123.4730  &  \\
\bottomrule
\end{tabular}
}
\end{table}

\subsection{Convergence w.r.t. discretization parameters}\label{sect:convg}
We study the relationship between the electrostatic solvation energy $E^{\rm s}$ and the discretization parameters, $\ell_{\rm max}$ and $N_{\rm leb}$.
For the sake of simplicity, we test a small molecule, namely formaldehyde with $4$ atoms (see Table \ref{tab:formaldehyde} for the geometry data).
First, we compute an ``exact'' electrostatic solvation energy with large discretization parameters $\ell_{\rm max}=25 $ and $N_{\rm leb} = 4334$.
This implies that we use $676$ basis functions and $4334$ integration points on each VDW sphere.

\begin{figure}
\centering
\includegraphics[width =5 in,clip,trim={0.in 0 1.2in 0}]{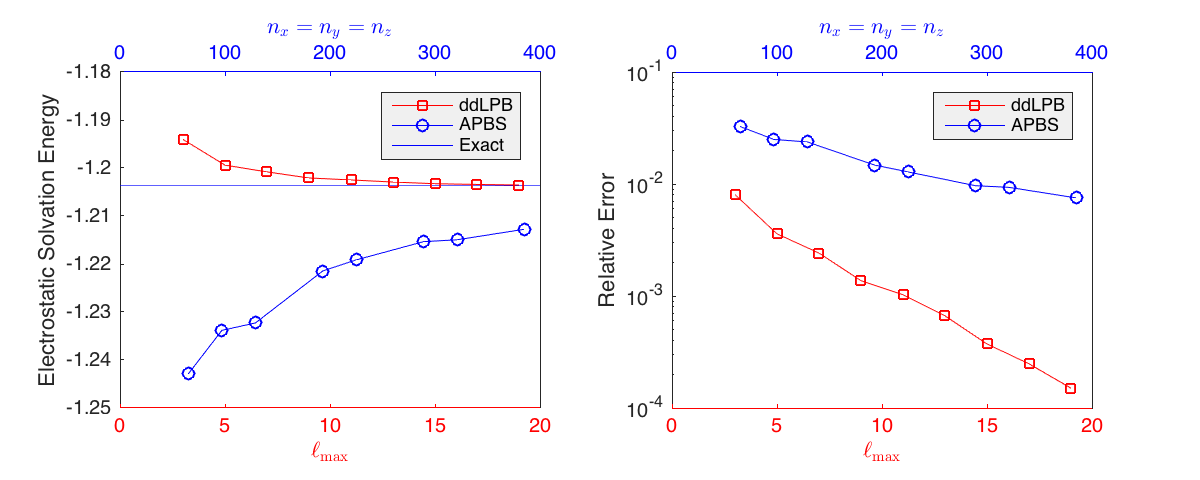}
\caption{The electrostatic solvation energies ({\rm kcal}/{\rm mol}, left) of formaldehyde and the relative errors (right). 
On the left-hand side, the blue line represents the ``exact'' electrostatic solvation energy; the red curve illustrates the energies computed by the ddLPB w.r.t. $\ell_{\rm max}$; the blue curve illustrates the energies computed by the APBS w.r.t. the number of grid points  $n_x=n_y=n_z$ in each axis direction.
}\label{fig:Test_lmax}
\end{figure}

The red curve in Figure \ref{fig:Test_lmax} illustrates how $E^{\rm s}$ computed by ddLPB varies w.r.t. the maximum degree of spherical harmonics $\ell_{\rm max}$, where $N_{\rm leb} = 4334$.
It is observed that the ddLPB provides systematically improvable approximations when $\ell_{\rm max}$ increases and we observed even exponential convergence of the energy w.r.t. $\ell_{\rm max}$. This allows to efficiently obtain an accuracy that is needed when the solvation model is coupled to a quantum mechanical description of the solute.
Further, we also run the APBS software for comparison, where the box size  is fixed to be $20\times 20\times 20$ (${\text \AA}^3$) and the grid dimension $n_x\times n_y\times n_z$ varies (with $n_x = n_y = n_z$ in $x$-, $y$-, $z$-axis). 
\orange In the input file, the molecular surface is chosen to be the VDW surface, by setting ${\tt srad = 0}$. 
We use the multi-grid solver by setting {\tt mg-manual}.
The other parameters are set as follows: {\tt gcent = mol 1, bcfl = mdh, chgm = spl4, sdens = 10, srfm = mol, swin = 0.3.}
\black
In Figure \ref{fig:Test_lmax}, the blue curve illustrates how $E^{\rm s}$ computed by the APBS varies w.r.t. $n_x$. 

\begin{figure}
\centering
\includegraphics[width =5 in,clip,trim={0.in 0 1.2in 0}]{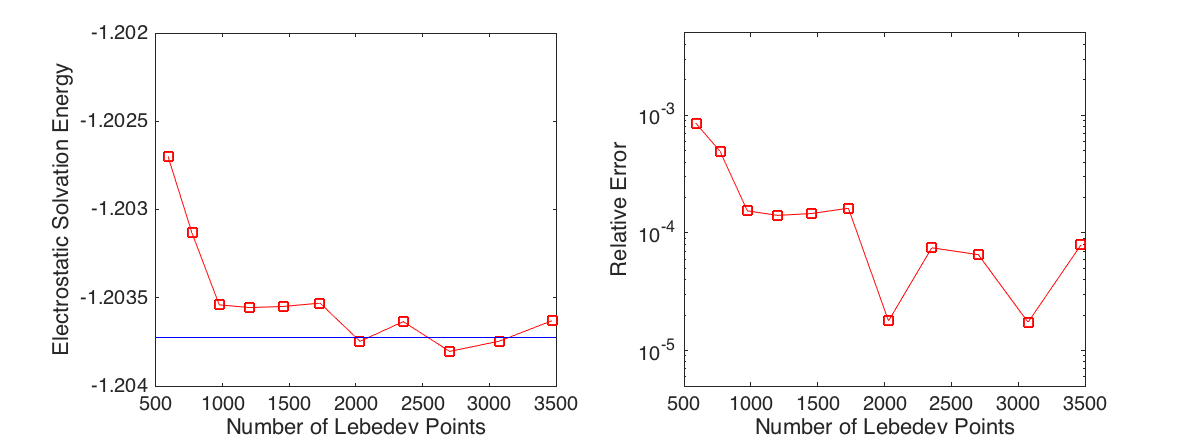}
\caption{ The electrostatic solvation energies ({\rm kcal}/{\rm mol}, left) of formaldehyde and the relative errors (right) w.r.t. $N_{\rm leb}$, computed by the ddLPB. 
On the left-hand side, the blue line represents the ``exact'' electrostatic solvation energy.
 }\label{fig:Test_Nleb}
\end{figure}

Figure \ref{fig:Test_Nleb} illustrates how $E^{\rm s}$ (computed by the ddLPB) varies w.r.t. the number of Lebedev quadrature points $N_{\rm leb}$, where $\ell_{\rm max} = 25$.
In fact, when $N_{\rm leb}$ is greater than $1000$, $E^{\rm s}$ varies very slightly (less than $0.04\%$), despite that it does not decay monotonically.

\begin{table}
\caption{Charges, centers $(x,y,z)$ and radii (\AA) of the $4$ atoms of formaldehyde.}\label{tab:formaldehyde}
\scriptsize
\centering
\scalebox{1.}{
\renewcommand\arraystretch{1.2}
\begin{tabular}{ccccc }
\toprule
Charge & $x$ & $y$ & $z$ & VDW-radius \\
\hline
0.08130 &  0.00000 &  0.00000 & -0.61750  & 2.11805 \\
\hline
-0.20542 &  0.00000 &  0.00000 &  0.75250  & 1.92500\\
\hline
0.06206 &  0.00000 &  0.93500 & -1.15750  & 1.58730 \\
\hline
0.06206 &  0.00000 & -0.93500 & -1.15750 &  1.58730 \\
\bottomrule
\end{tabular}
}
\end{table}

\black

\subsection{Varying the Debye-H\"uckel screening constants}
We now study the relationship between the electrostatic solvation energy and the Debye-H\"uckel screening constant $\kappa$.
\red 
On the continuous level, the solution of the Poisson-Boltzmann equation tends to the one of COSMO in the limit $\kappa\to \infty$.
Physically speaking, this is reasonable since the solvent becomes a perfect conductor as the ionic strength tends to $\infty$ and screens any charge from the solute.
On the other hand, the solution of PB equation tends to the solution of PCM in the limit $\kappa\to 0$.

\begin{figure}
\centering
\includegraphics[width =5 in,clip,trim={0.in 0 1.2in 0}]{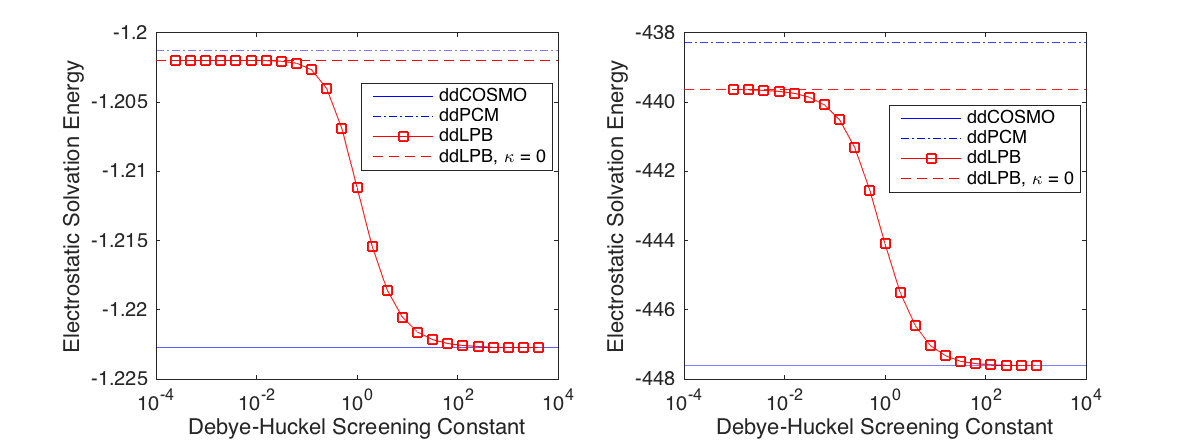}
\caption{The electrostatic solvation energy ({\rm kcal}/{\rm mol}) of formaldehyde (left) and 1etn (right) w.r.t. the Debye-H\"uckel screening constant $\kappa$. 
Solid blue line: ddCOSMO result; dashed blue line: ddPCM result; red curve: ddLPB result, $0< \kappa <\infty$; dashed red line: ddLPB result, $\kappa = 0$.
}\label{fig:Test_kappa}
\end{figure}

It is therefore natural to compare the ddLPB  with the ddCOSMO \cite{cances2013domain} and the ddPCM \cite{stamm2016new}.
In Figure \ref{fig:Test_kappa}, we plot the electrostatic solvation energies of formaldehyde (with $4$ atoms) and 1etn (PDB ID, with $141$ atoms) w.r.t. the Debye-H\"uckel screening constant.
The same discretization parameters $\ell_{\rm max} = 11$ and $N_{\rm leb} = 590$ are used for all three methods ddCOSMO, ddPCM and ddLPB.

For the limit $\kappa\to\infty$, the ddLPB result tends to the ddCOSMO result. 
This is consistent with the theory since the discretized equations for the ddLPB with $\kappa$ coincide with those for the ddCOSMO even after discretization. 
Indeed, the ddCOSMO method can been seen as a particular ddLPB method in the case of $\kappa = \infty$.

When $\kappa= 0$, we observe a difference between the ddLPB  and ddPCM results.
The reason is that the ddPCM discretizes the IEF-PCM \cite{cances1998new} directly whereas the ddLPB uses a discretization of a different integral formulation.
Since the continuous models are equivalent, the difference tends to zero when higher discretization parameters are used. 
\black

\subsection{Comparison with APBS}
In this part, we further compare the ddLPB solver with the widely-used software, APBS.
\orange 
The same parameters are used as in Section \ref{sect:convg}, except the grid spacings and dimensions.
\black
The test is performed on a MacBook Pro with a 2.5 GHz Intel Core i7 processor and we consider the protein 1etn with $141$ atoms as test case.

\begin{figure}
\centering
\includegraphics[width =5 in,clip,trim={0in 0 1.2in 0}]{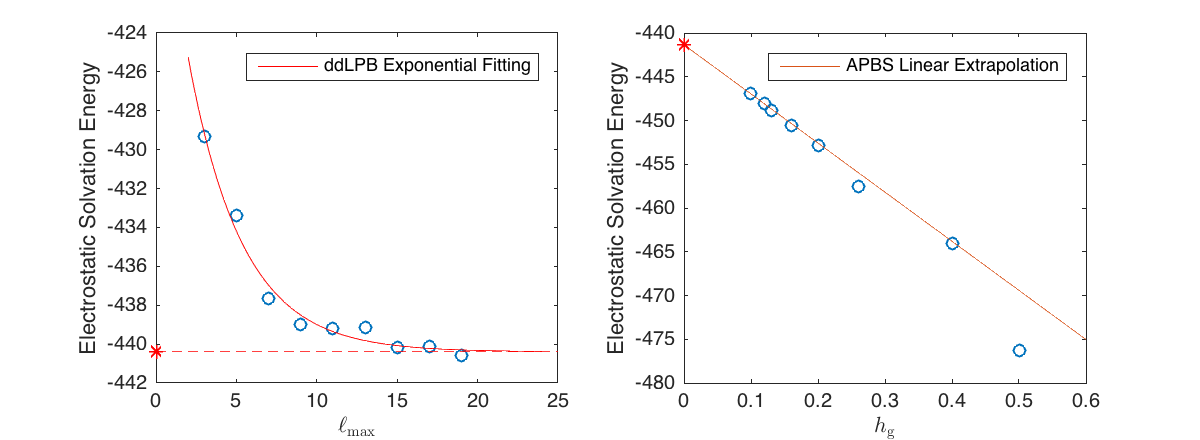}
\caption{\red Estimation of the electrostatic solvation energy of 1etn (141 atoms), based on the exponential data fitting for the ddLPB results (left) and the linear extrapolation for the APBS results (right) in Table \ref{tab:ddLPBvsAPBS}. 
Here, $h_{\rm g}$ represents the grid spacing in the APBS. 
On the left, the red curve plots the fitting function $y = -440.3956 + 27.4584 ~{\rm exp}({-0.2974x})$, which tends to $E^*_{ddlpb} = -440.3956$ at the infinity.
The dashed horizontal line is the asymptote of the curve.
On the right, the red line crosses the leftmost two energies and intersects the $y$-axis at the star maker $E^*_{apbs} = -441.3488$.
}\label{fig:extrapolation}
\end{figure}

Table \ref{tab:ddLPBvsAPBS} illustrates the ddLPB and APBS results for different discretization parameters, including the electrostatic solvation energy, the relative error, the number of iterations, the run time and the memory.
\red
To compute the relative error, two ``exact'' electrostatic solvation energies are computed respectively using an exponential fitting for ddLPB and a linear extrapolation for APBS, as illustrated in Figure \ref{fig:extrapolation}. 

In fact, we have observed an exponential convergence of the energy w.r.t. $\ell_{\rm max}$ for formaldehyde in Section \ref{sect:convg}.
It appears therefore consistent to apply an exponential data fitting for the ddLPB-energies.
More precisely, we use the function {\tt fit} in Matlab where the fit type is set to $ y=a+b\cdot {\rm exp}(-c\cdot x)$. 
The ``exact'' electrostatic solvation energy $E^{*}_{\rm ddlpb}$ is obtained as the coefficient $a$ when the fitting function is figured out.
For the APBS, we use the linear extrapolation procedure introduced in \cite{geng2013treecode}. 
We first plot the APBS energies w.r.t. $h_{\rm g}$ and then draw a line crossing the leftmost two energies at $h_{\rm g} = 0.1$ and $h_{\rm g}  = 0.12$.
As a consequence, as $h_{\rm g}$ tends to zero, this line crosses the $y$-axis at an ``exact'' electrostatic solvation energy $E^{*}_{\rm apbs}$.
\black

In Table \ref{tab:ddLPBvsAPBS}, we can observe that the ddLPB usually cost less memory than the APBS due to the nuclear-centered spectral-type basis functions similar to atomic orbitals. 
Furthermore, from the relative errors obtained by extrapolation, one observes that the ddLPB results are more accurate in this example, as also observed in Figure~\ref{fig:Test_lmax} which can be explained that APBS is a first order method while ddLPB shows exponential convergence for the energy.
\black

\begin{table}
\caption{ddLPB (top) and APBS results (bottom) for the protein 1etn. $N_{\rm iter}$ represents the number of outer iterations in the ddLPB. $h_{\rm g}$ and $N_{\rm g}$ represent the grid spacing and the grid dimension in the APBS. {\red $E^{*}_{\rm ddlpb}$ and $E^{*}_{\rm apbs}$ represent the  ``exact'' electrostatic solvation energies computed respectively from an exponential fitting for ddLPB and a linear extrapolation for APBS, as illustrated in Fig. \ref{fig:extrapolation}. }
}\label{tab:ddLPBvsAPBS}
\scriptsize
\centering
\scalebox{1.}{
\renewcommand\arraystretch{1.2}
\begin{tabular}{c cccccccc}
\toprule
{\bfseries $\ell_{\rm max}$}& {\bfseries $N_{\rm leb}$} & $E^{\rm s} ({\rm kcal}/{\rm mol})$ &{\tt RE} & $N_{\rm iter}$ & Run time ({\rm s}) & {Memory} ({\rm MB})  \\
\midrule
3 & 26 & -429.3457 & 2.5091e-02 &  9 & 1 & 7  \\
\hline
5 & 50 & -433.3903& 1.5907e-02 &  11  & 2& 30  \\
\hline
7 &86 & -437.6449 &  6.2460e-03 &  11  & 8 & 91  \\
\hline
9 & 146 & -438.9725 & 3.2314e-03 &  11  & 22 & 241  \\
\hline
11 & 194&  -439.1880 & 2.7421e-03 &  14 & 63 & 461  \\
\hline
13 & 266&  -439.1629 & 2.7991e-03 &  18 &  198 &861 \\
\hline
15 & 350&  -440.1754 & 5.0001e-04 &  12 & 206 &893  \\
\hline
17 & 434&  -440.1017 & 6.7735e-04 &  13 & 362 &1398  \\
\hline
 \multicolumn{2}{c}{$E^{*}_{\rm ddlpb}$} &  -440.3956  & &   &  &    \\
\bottomrule
\toprule
{\bfseries $h_{\rm g}$}& {\bfseries $N_{\rm g}$} & $E^{\rm s} ({\rm kcal}/{\rm mol})$ &{\tt RE} & $N_{\rm iter}$ & Run time ({\rm s}) & {Memory} ({\rm MB})  \\
\midrule 
0.5 & $65^3$ & -476.3143  &  7.9224e-02 & -- &  1 &63    \\
\hline
0.4 & $97^3$  & -464.0155  &  5.1358e-02 & -- &   3 &203   \\
\hline
0.26 & $129^3$ &  -457.4919  &  3.6577e-02 & --  &  8  &475   \\
\hline
0.2 & $193^3$ &  -452.8998  &  2.6172e-02 & -- &  29  &1606   \\
\hline
0.16 & $225^3$ &  -450.5827  &  2.0922e-02 & -- & 46   &2572   \\
\hline
0.13 & $289^3$ &  -448.8017  &  1.6887e-02 & -- &   113 &5608 \\
\hline
0.12 & $321^3$ &  -448.0838 &  1.5260e-02 & -- &  190  &7819  \\
\hline
0.1 & $385^3$ &  -446.9613 &  1.2717e-02 & -- &  525  &14044\\
\hline
 \multicolumn{2}{c}{$E^{*}_{\rm apbs}$} &  -441.3488  & &   &  &    \\
\bottomrule
\end{tabular}
}
\end{table}

\subsection{Computational cost}
To study the computational cost of the ddLPB, we test a set of $24$ proteins with the following PDB IDs: 1ajj, 1ptq, 1vjw, 1bor, 1fxd, 1sh1, 1hpt, 1fca, 1bpi, 1r69, 1bbl, 1vii, 2erl, 451c, 2pde, 1cbn, 1frd, 1uxc, 1mbg, 1neq, 1a2s, 1svr, 1o7b, 1a63.
The discretization parameters of the ddLPB are set to $\ell_{\rm max} = 5$ and $N_{\rm leb} = 50$.
For the comparison, we also run the TABI-PB solver \cite{geng2013treecode}, which however works on the SES, not the VDW surface.
To have approximately the same degree of freedom, we set the density of vertices in the TABI-PB to be $5$ per ${\text \AA}^{2}$.
In the TABI-PB, the order of Taylor approximation is set to $1$, the MAC parameter is set to $0.4$ and the GMRES tolerance (relative residual error) is set to $10^{-4}$.
\red
Both the TABI-PB and the ddLPB use the default parameters given at the beginning of Section \ref{sect:numerical}, that is, $\varepsilon_1=1,~\varepsilon_2 = 78.54$ and $\kappa = 0.1040~ {\text \AA}^{-1}$.
\black

\begin{figure}
\centering
\includegraphics[width =5 in,clip,trim={0.in 0 1.2in 0}]{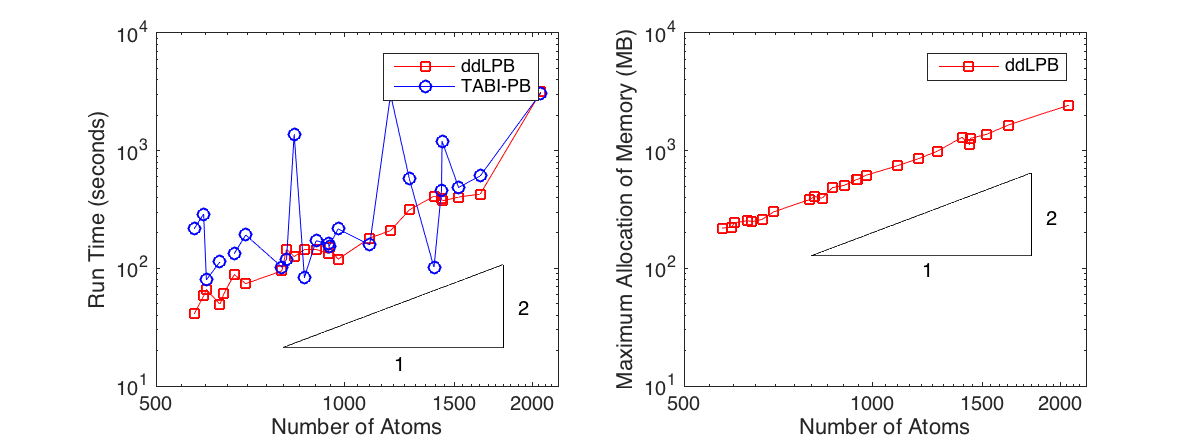}
\caption{Run time (left) and maximum allocation of memory (right) of the ddLPB and the TABI-PB w.r.t. the number of atoms, for different proteins in the test set. The detailed memory information of the TABI-PB is not available and the memory is finally released.
}\label{fig:cost}
\end{figure}

Figure \ref{fig:cost} illustrates the run time and the maximum allocation of memory of the ddLPB and the TABI-PB for different proteins in the test set.
Table \ref{tab:info} provides more information including the degree of freedom and the number of iterations. As a general observation, we see that the number of iteration to reach the tolerance is more stable in ddLPB.
Notice that for the protein 1cbn, the MSMS in the TABI-PB fails to generate a suitable mesh and the TABI-PB stops without returning the energy.
Further, for the protein 1bor and 1neq, the GMRES in the TABI-PB reaches the maximum number of iterations 100, before reaching the residual tolerance $10^{-4}$.
The ddLPB reaches the tolerance ${\tt Tol} = 10^{-4}$ within $20$ outer iterations for all proteins in the test set.

%

\begin{table}
\caption{Details on the ddLPB and the TABI-PB for the test set of proteins, including the degree of freedom and the number of iterations. {\orange $M$ denotes the number of atoms}. }\label{tab:info}
\scriptsize
\centering
\scalebox{0.89}{
\renewcommand\arraystretch{1.2}
\begin{tabular}{cccccc}
\toprule
 & &  \multicolumn{2}{c}{Degree of freedom}  & \multicolumn{2}{c}{Iteration} \\ \cline{3-6}
PDB & $M$  & ddLPB & TABI  & ddLPB & TABI \\
\midrule
1ajj & 602& 43344 & 41420 &  15 & 10 \\
\hline
1ptq & 795& 57240 & 56898 &  12 & 9 \\ 
\hline
1vjw & 946& 68112 & 56084 & 12  & 13 \\ 
\hline
1bor &832 & 59904 & 58076 & 16 & 100 \\ 
\hline
1fxd & 978& 70416 & 62844 & 9 & 14 \\
\hline
1sh1 & 696& 50112 & 55604 & 14 & 15 \\ 
\hline
1hpt & 945 & 68040 & 65572 &  12 & 10\\ 
\hline
1fca & 864 & 62208 & 49744 & 14 & 9 \\ 
\hline
 1bpi & 1393 & 103320 & 64482 &13  & 7\\ 
\hline
1r69 & 1099 & 79128 & 62648 & 11 & 11 \\ 
\hline
1bbl & 576 & 41472 & 53654 & 10 & 17 \\ 
\hline
1vii & 596 & 42912 & 50874 & 14 & 24 \\ 
\bottomrule
\end{tabular}
\begin{tabular}{cccccc}
\toprule
 &  & \multicolumn{2}{c}{Degree of freedom}  & \multicolumn{2}{c}{Iteration}  \\ \cline{3-6}
PDB & $M$ & ddLPB & TABI  & ddLPB  & TABI  \\
\midrule
2erl & 633 & 45576 & 46642 & 10 & 11 \\ 
\hline
451c & 1435 & 103320 & 89286 & 13 & 45 \\ 
\hline
2pde & 667 & 48024 & 48430 & 18 & 12 \\ 
\hline
1cbn & 642 & 46224 &  -- & 12  & --\\ 
\hline
1frd & 1652 & 118944 & 90772 & 11  & 25 \\
\hline
1uxc & 809 & 58248 & 58164 & 19 & 10 \\ 
\hline
1mbg & 902 & 64944 & 64190 & 14 & 12 \\
\hline
1neq & 1187 & 85464 & 98964 & 11 & 100 \\ 
\hline
1a2s & 1272 & 91584 & 93574 & 15 & 22 \\
\hline
 1svr & 1432 & 103104 & 100546 & 15 & 16 \\ 
\hline
1o7b & 1525 & 109800 & 103920 & 13 & 17 \\ 
\hline
1a63 & 2065 & 148680 & 144796 & 13 & 74 \\ 
\bottomrule
\end{tabular}
}
\end{table}

\begin{remark}
At this moment, we haven't yet employed acceleration techniques in the ddLPB implementation, while the TABI-PB features the ``treecode'' acceleration technique \cite{geng2013treecode}.
\end{remark}

\subsection{Graphical illustration}
\begin{figure}
\centering
\includegraphics[width =1.66 in,clip,trim={1.7in 1in 1in 0.5in}]{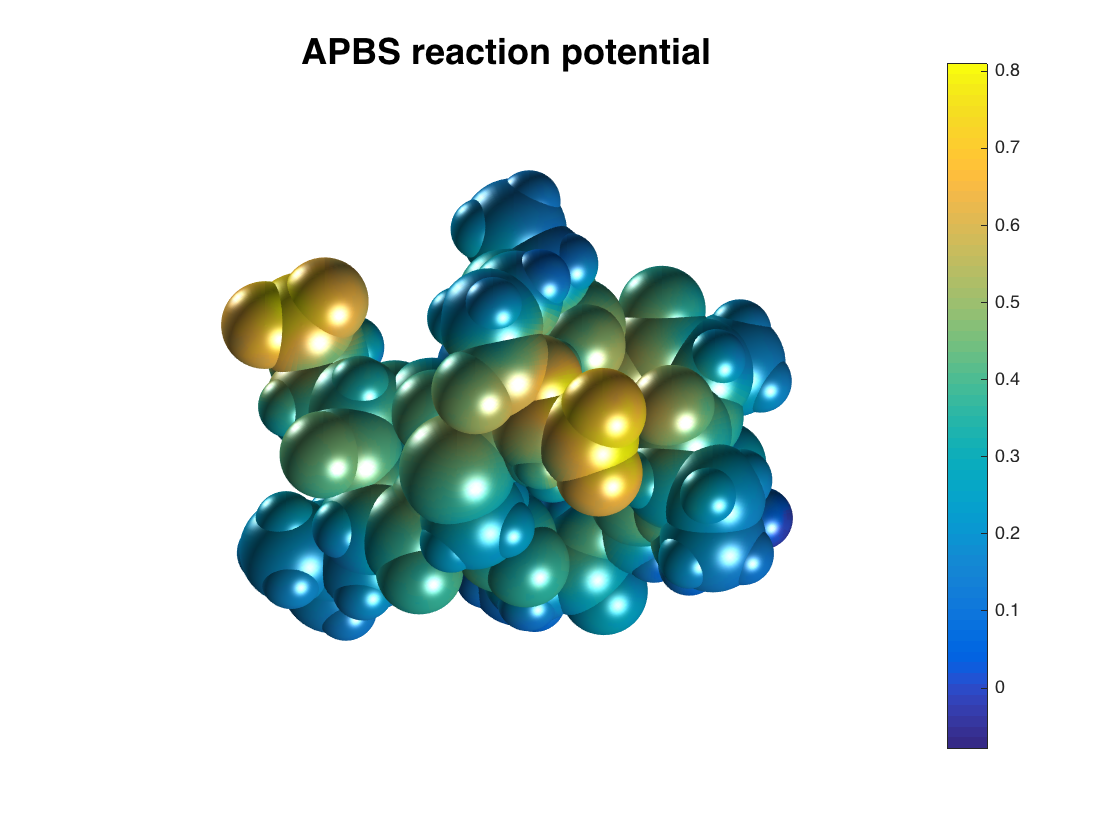}
\includegraphics[width =1.66 in,clip,trim={1.7in 1in 1in 0.5in}]{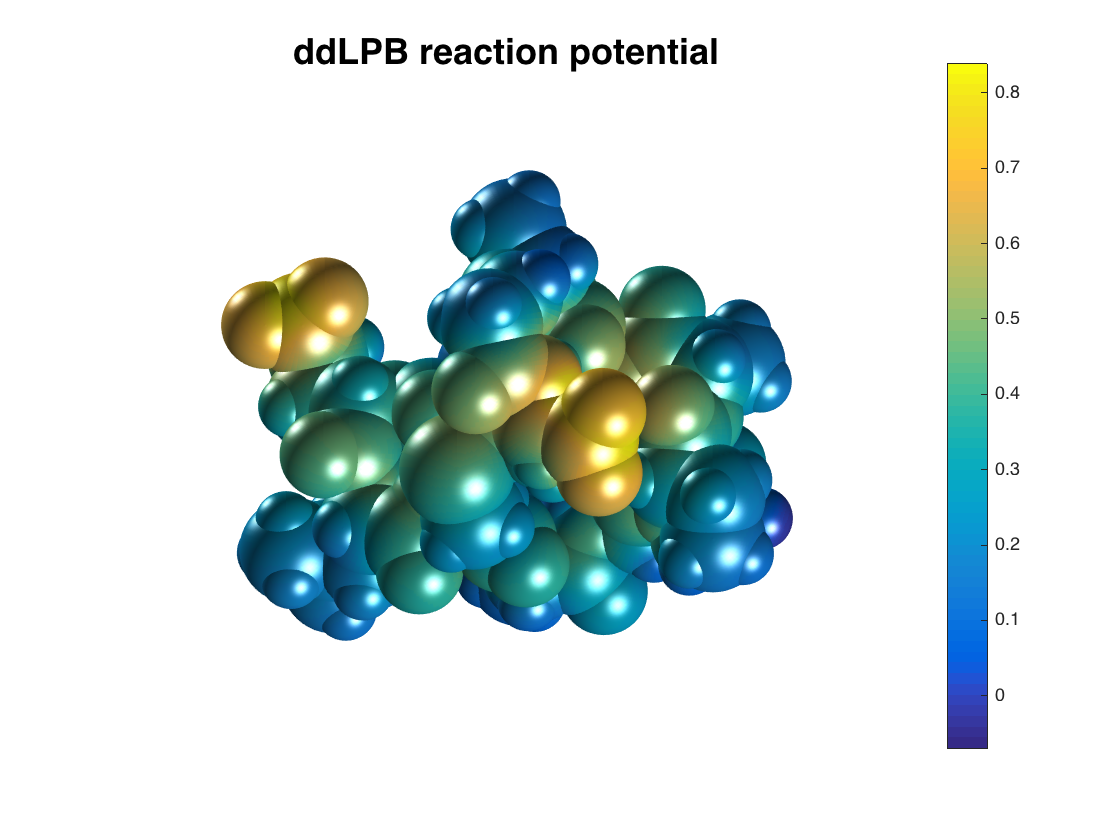}
\includegraphics[width =1.66 in,clip,trim={1.4in 1in 1in 0.5in}]{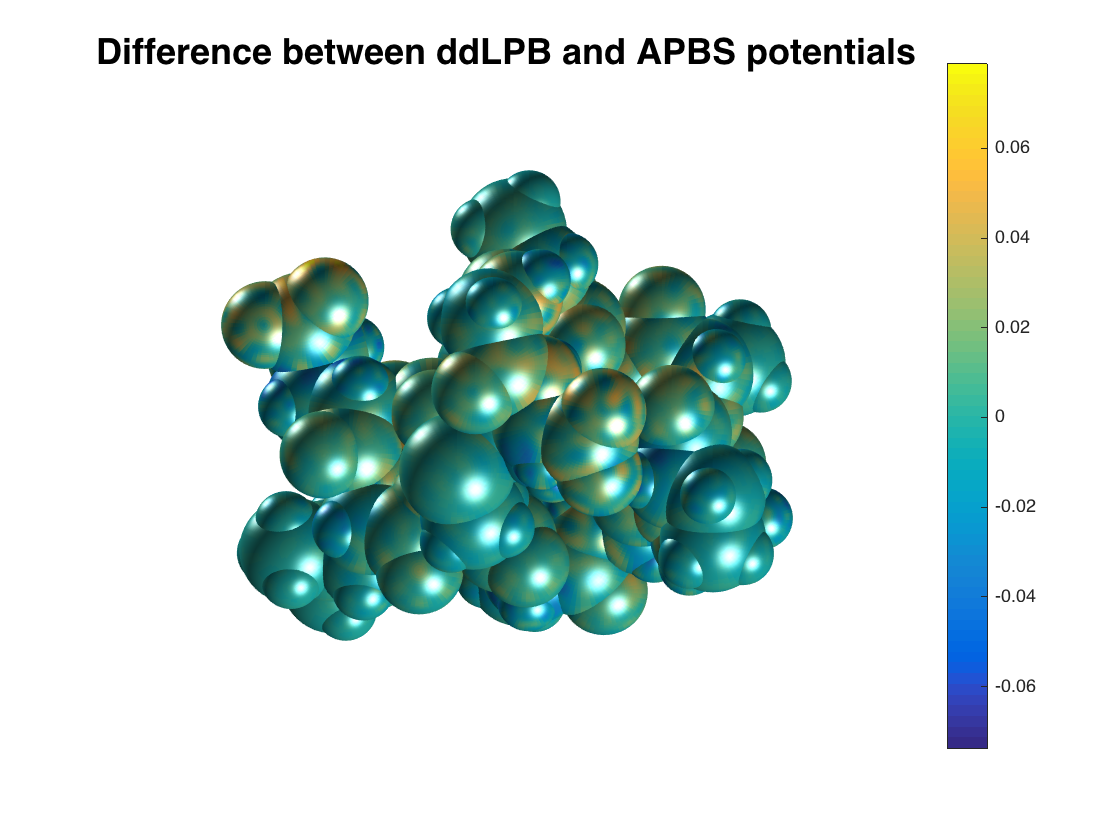}
\caption{\blue Reaction potentials ({\mbox e}/{\mbox \AA} = 561.91 {\mbox kT}/{\mbox e})  on the VDW surface of 1etn with 141 atoms: APBS result (left), ddLPB result (middle) and difference (right). The following parameters are used: grid dimension $193\times 193\times 193$ and grid spacing $0.2\times 0.2\times 0.2$ for the APBS, $\ell_{\rm max} = 11$ and $N_{\rm leb} = 194$ for the ddLPB. The other parameters are the default ones as used previously. }\label{fig:graph1}
\end{figure}

\begin{figure}
\centering
\includegraphics[width = 1.45in,clip,trim={3.5in 0.8in 1.in 0.4in}]{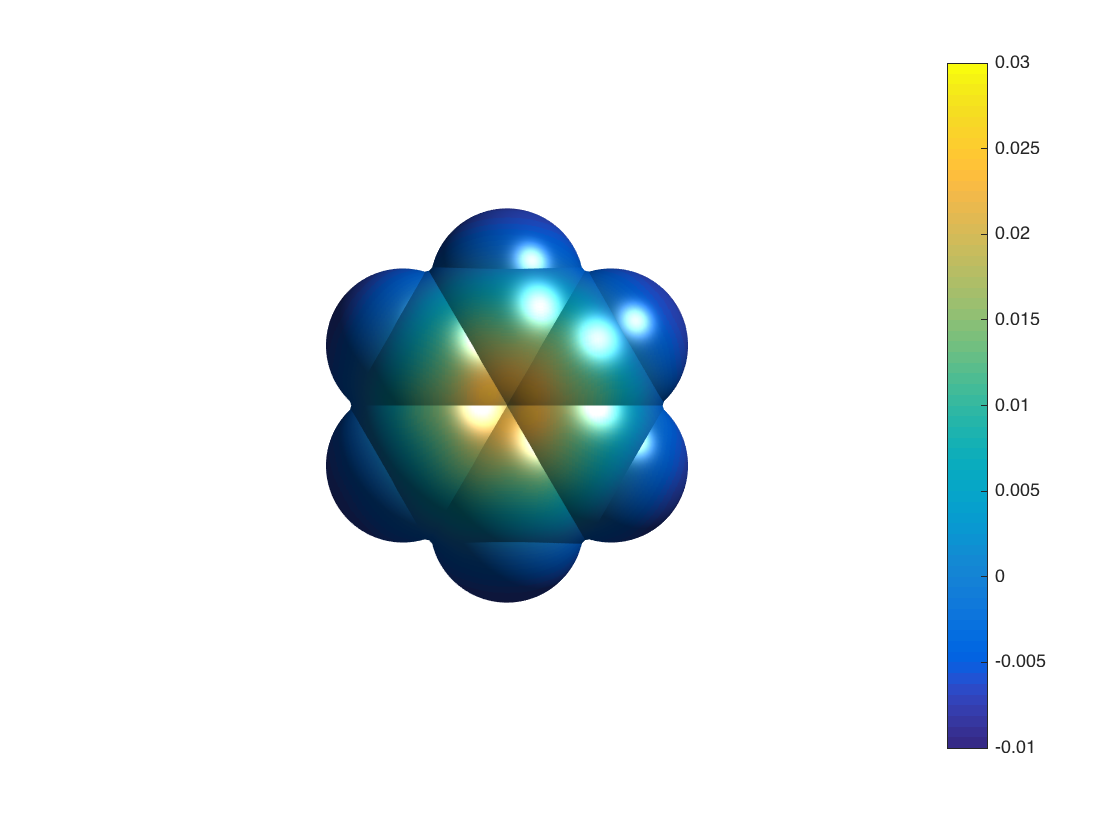}
\includegraphics[width = 1.56in,clip,trim={3.0in 0.8in 1.in 0.5in}]{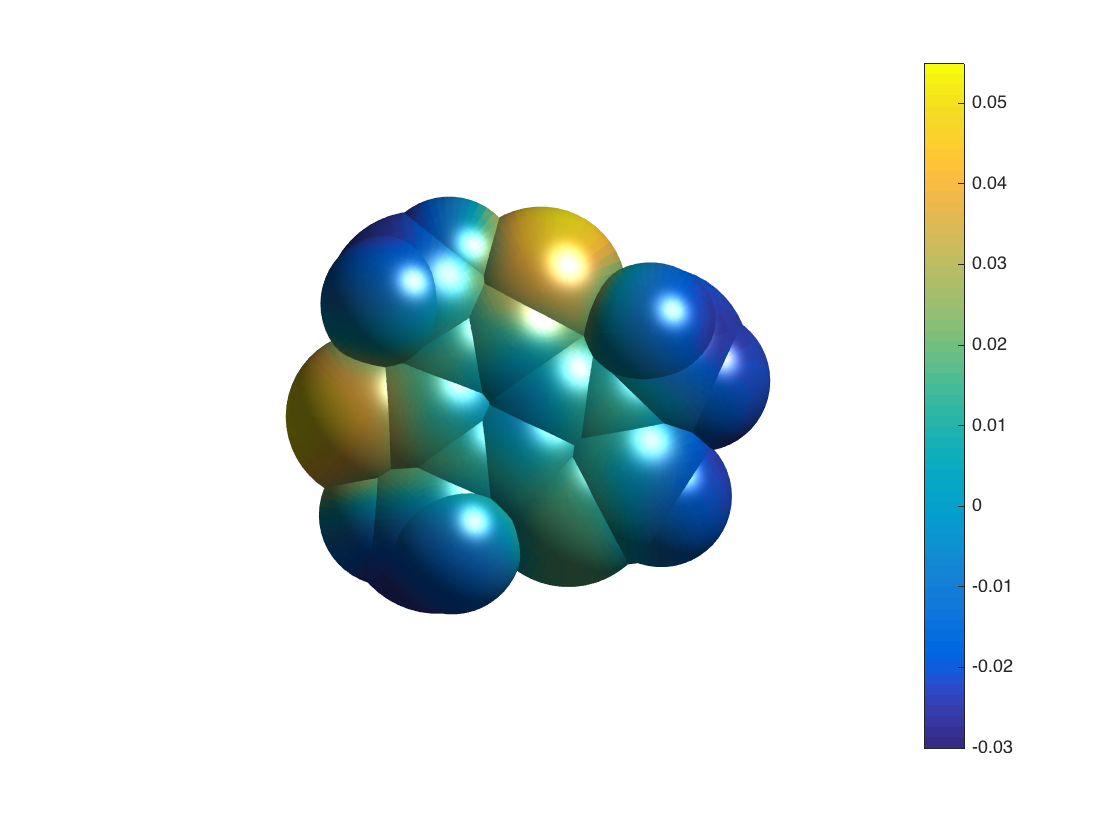}
\includegraphics[width = 1.78in,clip,trim={1.5in 1in 1.in 0.6in}]{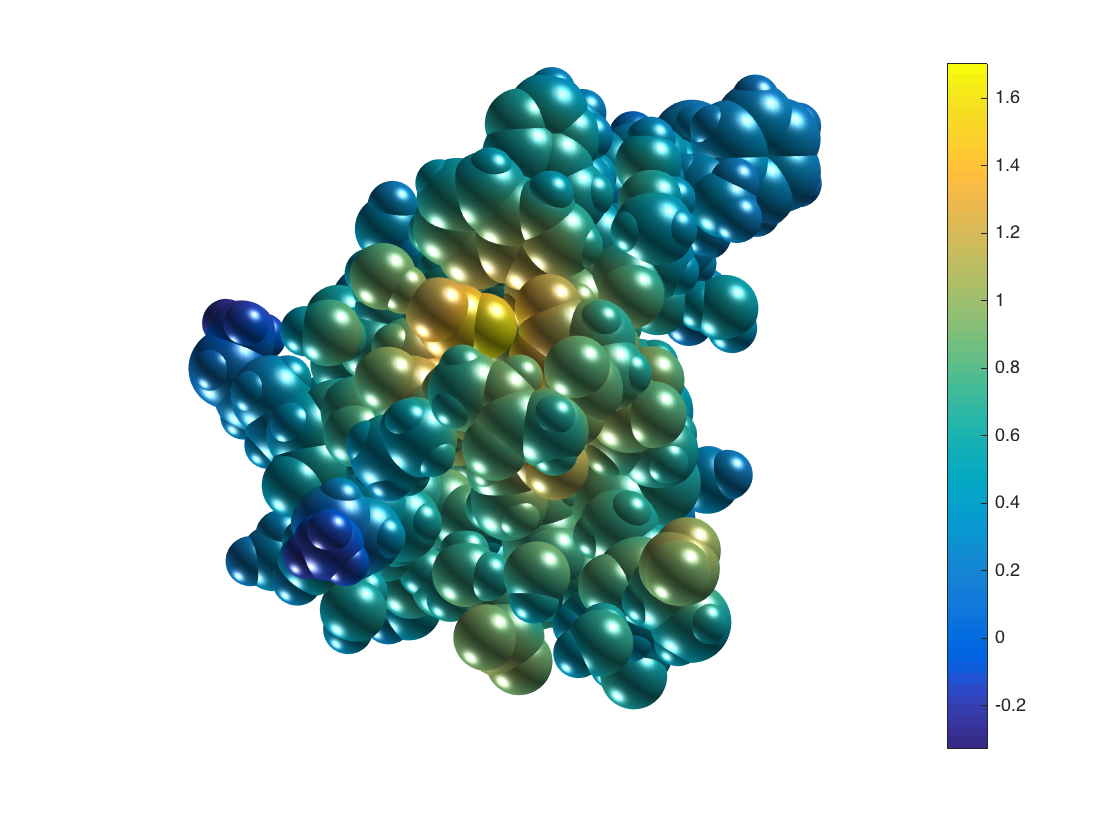}
\caption{\blue Reaction potential ({\mbox e}/{\mbox \AA}) on the VDW surfaces of benzene (left), caffeine (middle) and 1ajj (right) respectively with $4$, $24$ and $602$ atoms. The following parameters are used: $\ell_{\rm max} = 11$ and $N_{\rm leb} = 194$ for benzene and caffeine, $\ell_{\rm max} = 7$ and $N_{\rm leb} = 86$ for the protein 1ajj. The other parameters are the default ones as used previously.}\label{fig:graph2}
\end{figure}

Finally, we provide some graphical illustrations of the reaction potential $\psi_{\rm r}$ on the VDW surface. 
\blue In Figure \ref{fig:graph1}, we illustrate the reaction potentials of 1etn computed by the APBS and the ddLPB, and their potential difference on the VDW surface. 
It is observed that the difference is usually small over the surface. 
In Figure \ref{fig:graph2}, we present the reaction potentials of two small molecules benzene and caffeine, and the protein 1ajj.
We notice that the reaction potential of benzene has rotational symmetry, which matches its geometrical structure.
\black

\section{Conclusion}
In this paper, we proposed a domain decomposition method for the Poisson-Boltzmann solvation model that shows exponential convergence in the energy w.r.t to the number of basis functions employed. This allows to reach a precision which enable this method to couple it with models on the level of quantum mechanics.
\black 

The original problem defined in $\mathbb R^3$ is first transformed into two coupled equations defined in the bounded solute cavity, based on potential theory.
Then, the Schwarz domain decomposition method was used to solve these two problems by decomposing the solute cavity into balls.
In consequence, we developed two direct single-domain solvers respectively for solving the Laplace equation and the HSP equation defined in the unit ball, which becomes easy to tackle by using the spherical harmonics in the angular direction.
Taking into account the coupling conditions allowed then to obtain a global linear system. 
A series of numerical results have been presented to show the performance of the ddLPB method.

In the future, we will focus on accelerating the ddLPB to make it suitable to very large molecules based on linear scaling acceleration techniques such as the Fast Multipole Method (FMM).
In addition, we will embed the ddLPB solver in software packages which simulate the computation of the solute on the level of theory of quantum mechanics or molecular dynamics.

\section{Acknowledgements}
We would like to thank F. Lipparini for discussion and guidance with the implementation in Fortran. 
We also thank the other members in our DD-family \cite{ddfamily} for the usual fruitful discussions, including E. Canc\`es, L. Lagard\`ere, J.P. Piquemal and B. Mennucci.

\bibliographystyle{unsrt}
\bibliography{bibfile.bib}

\appendix
\section{Computation of $\mathbf C_1,~\mathbf C_2$ and $F_0$} \label{append:1}
For each $\Gamma_i^{\rm e}$, we first define a square matrix $\mathbf P_{\chi^{\rm e}_i}$ of dimension $(\ell_{\rm max}+1)^2\times(\ell_{\rm max}+1)^2$ for each $\chi^{\rm e}_i$, the $(\ell_0 m_0 , \ell' m')$-th element of which is defined as
\begin{equation}\label{eq:P}
[\mathbf P_{\chi^{\rm e}_i}]_{\ell_0 m_0}^{\ell'm'} \coloneqq \sum_{n = 1}^{N_{\rm leb}} w_{n}^{\rm leb}\,\chi_i^{\rm e}(\mathbf x_i+r_i \mathbf s_n)\, Y_{\ell_0}^{m_0}(\mathbf s_n)\,Y_{\ell'}^{m'}(\mathbf s_n),
\end{equation}
where $0\leq \ell_0\leq\ell_{\rm max},-\ell_0\leq m_0\leq\ell_0, 0\leq \ell'\leq\ell_{\rm max},-\ell'\leq m'\leq\ell'.$
Based on Eq. \eqref{eq:partial_psi_r}, we can approximate $\chi_i^{\rm e}\partial_{\mathbf n}\psi_{\rm r}$ (defined on $\Gamma_i$) by a linear combination of spherical harmonics $\{ Y_{\ell_0}^{m_0} \}$ with $0\leq \ell_0\leq\ell_{\rm max},-\ell_0\leq m_0\leq\ell_0$ as follows
\begin{equation}
\chi_i^{\rm e}\partial_{\mathbf n}\psi_{\rm r}(r_i,\mathbf s) = \sum_{\ell_0 = 0}^{\ell_{\rm max}} \sum_{m_0=-\ell_0}^{\ell_0} c_{{\rm r},\ell_0 m_0}\, Y_{\ell_0}^{m_0}(\mathbf s),\quad \mathbf s\in \mathbb S^2,
\end{equation}
where the coefficient $c_{{\rm r},\ell_0 m_0}$ is computed by the Lebedev quadrature rule as follows
\begin{equation}
\begin{array}{r@{}l}
\begin{aligned}
c_{{\rm r},\ell_0 m_0} = \sum_{\ell'=0}^{\ell_{\rm max}}\sum_{m'=-\ell'}^{\ell'} \, [\mathbf P_{\chi^{\rm e}_i}]_{\ell_0 m_0}^{\ell'm'} \, \frac{\ell'}{r_i} \,  [X_{\rm r}]_{i \ell' m'}.
\end{aligned}
\end{array}
\end{equation}
\begin{remark}
By writing $\chi_i^{\rm e}\partial_{\mathbf n}\psi_{\rm r}$ as a linear combination of spherical harmonics, the single-layer potential $\widetilde{\mathcal S}_{\kappa,\Gamma_i}$ can act on it conveniently.
\end{remark}
For an arbitrary Lebedev point $\mathbf x_j + r_j\mathbf s_n =  \mathbf x_i + r_{ijn}\mathbf s_{ijn} \in \Gamma_j^{\rm e}$, we can then compute as follows
\begin{equation}
\begin{array}{r@{}l}
\begin{aligned}
&\left( \widetilde{\mathcal S}_{\kappa,\Gamma_i}\chi_i^{\rm e}\partial_{\mathbf n}\psi_{{\rm r}} \right)(\mathbf x_j + r_j\mathbf s_n) \\
&= \sum_{\ell_0 = 0}^{\ell_{\rm max}} \sum_{m_0=-\ell_0}^{\ell_0} c_{{\rm r},\ell_0 m_0}\, \left( \widetilde{\mathcal S}_{\kappa,\Gamma_i}Y_{\ell_0}^{m_0}\right)(\mathbf x_j + r_j\mathbf s_n)\\
& =  \sum_{\ell_0 = 0}^{\ell_{\rm max}} \sum_{m_0=-\ell_0}^{\ell_0} c_{{\rm r},\ell_0 m_0}\,  \left( \frac{ i'_{\ell_0}\left({ r_i}\right)}{i_{\ell_0}\left({}{r_i}\right)} -  \frac{k'_{\ell_0}\left({r_i}\right)}{k_{\ell_0}\left({}{r_i}\right)} \right)^{-1} \,\frac{ k_{\ell_0}\left({ r_{ijn}}\right)}{k_{\ell_0}\left({}{r_i}\right)}\, Y_{\ell_0}^{m_0}(\mathbf s_{ijn}) \\
& = \sum_{\ell'=0}^{\ell_{\rm max}}\sum_{m'=-\ell'}^{\ell'} [\mathbf Q]_{i\ell'm'}^{jn} \, \frac{\ell'}{r_i} \,  [X_{\rm r}]_{i \ell' m'},
\end{aligned}
\end{array}
\end{equation} 
where $\mathbf Q$ is a matrix of dimension $M(\ell_{\rm max}+1)^2\times M N_{\rm leb}$, with the $(i\ell'm',jn)$-th element defined by
\begin{equation}\label{eq:Q}
[\mathbf Q]_{i\ell'm'}^{jn} \coloneqq \sum_{\ell_0 = 0}^{\ell_{\rm max}} \sum_{m_0=-\ell_0}^{\ell_0} [\mathbf P_{\chi^{\rm e}_i}]_{\ell_0 m_0}^{\ell'm'}\,  \left( \frac{i'_{\ell_0}\left({ r_i}\right)}{i_{\ell_0}\left({r_i}\right)} -  \frac{k'_{\ell_0}\left({ r_i}\right)}{k_{\ell_0}\left({}{r_i}\right)} \right)^{-1} \,\frac{ k_{\ell_0}\left({ r_{ijn}}\right)}{k_{\ell_0}\left({}{r_i}\right)}\, Y_{\ell_0}^{m_0}(\mathbf s_{ijn}).
\end{equation}
Therefore, we have the $(j\ell m,i\ell'm')$-th element of $\mathbf C_1$ as follows
\begin{equation}\label{eq:C_1}
[\mathbf C_1]_{j\ell m}^{i\ell'm'} = \frac{\varepsilon_1}{\varepsilon_2}\left(\sum_{n = 1}^{N_{\rm leb}} w_n^{\rm leb} \chi_j^{\rm e}(\mathbf x_j +r_j\mathbf s_n) Y_{\ell}^m(\mathbf s_n) [\mathbf Q]_{i\ell'm'}^{jn} \, \frac{\ell'}{r_i}\right).
\end{equation}

Similarly, based on Eq. \eqref{eq:partial_psi_e}, we can approximate $\chi_i^{\rm e}\partial_{\mathbf n}\psi_{\rm e}$ (defined on $\Gamma_i$) by another linear combination of spherical harmonics $\{ Y_{\ell_0}^{m_0} \}$ with $0\leq \ell_0\leq\ell_{\rm max},-\ell_0\leq m_0\leq\ell_0$ as follows
\begin{equation}
\chi_i^{\rm e}\partial_{\mathbf n}\psi_{\rm e}(r_i,\mathbf s) = \sum_{\ell_0 = 0}^{\ell_{\rm max}} \sum_{m_0=-\ell_0}^{\ell_0} c_{{\rm e},\ell_0 m_0}\, Y_{\ell_0}^{m_0}(\mathbf s),\quad \mathbf s\in \mathbb S^2,
\end{equation}
where
\begin{equation}
\begin{array}{r@{}l}
\begin{aligned}
c_{{\rm e},\ell_0 m_0} = \sum_{\ell'=0}^{\ell_{\rm max}}\sum_{m'=-\ell'}^{\ell'} [\mathbf P_{\chi^{\rm e}_i}]_{\ell_0 m_0}^{\ell'm'} \,\frac{i'_{\ell'}\left({ r_i}\right)}{i_{\ell'}\left({}{r_i}\right)} \, [X_{\rm e}]_{i \ell' m'}.
\end{aligned}
\end{array}
\end{equation}
For an arbitrary Lebedev point $\mathbf x_j + r_j\mathbf s_n =  \mathbf x_i + r_{ijn}\mathbf s_{ijn} \in \Gamma_j^{\rm e}$, we can then compute
\begin{equation}
\begin{array}{r@{}l}
\begin{aligned}
&\left( \widetilde{\mathcal S}_{\kappa,\Gamma_i}\chi_i^{\rm e}\partial_{\mathbf n}\psi_{{\rm e}} \right)(\mathbf x_j + r_j\mathbf s_n) \\
&= \sum_{\ell_0 = 0}^{\ell_{\rm max}} \sum_{m_0=-\ell_0}^{\ell_0} c_{{\rm e},\ell_0 m_0}\, \left( \widetilde{\mathcal S}_{\kappa,\Gamma_i}Y_{\ell_0}^{m_0}\right)(\mathbf x_j + r_j\mathbf s_n)\\
& =  \sum_{\ell_0 = 0}^{\ell_{\rm max}} \sum_{m_0=-\ell_0}^{\ell_0} c_{{\rm e},\ell_0 m_0}\, \left( \frac{i'_{\ell_0}\left({ r_i}\right)}{i_{\ell_0}\left({}{r_i}\right)} -  \frac{ k'_{\ell_0}\left({ r_i}\right)}{k_{\ell_0}\left({}{r_i}\right)} \right)^{-1} \,\frac{ k_{\ell_0}\left({ r_{ijn}}\right)}{k_{\ell_0}\left({}{r_i}\right)}\, Y_{\ell_0}^{m_0}(\mathbf s_{ijn}) \\
& = \sum_{\ell'=0}^{\ell_{\rm max}}\sum_{m'=-\ell'}^{\ell'} [\mathbf Q]_{i\ell'm'}^{jn} \, \frac{ i'_{\ell'}\left({ r_i}\right)}{i_{\ell'}\left({}{r_i}\right)} \, [X_{\rm e}]_{i \ell' m'}.
\end{aligned}
\end{array}
\end{equation} 
This yields that
\begin{equation}\label{eq:C_2}
[\mathbf C_2]_{j\ell m}^{i\ell'm'} = - \left(\sum_{n = 1}^{N_{\rm leb}} w_n^{\rm leb} \chi_j^{\rm e}(\mathbf x_j +r_j\mathbf s_n) Y_{\ell}^m(\mathbf s_n) [\mathbf Q]_{i\ell'm'}^{jn} \, \frac{ i'_{\ell'}\left({ r_i}\right)}{i_{\ell'}\left({}{r_i}\right)} \right).
\end{equation}

In addition, since $\partial_{\mathbf n}\psi_0$ is known, for an arbitrary Lebedev point $\mathbf x_j + r_j\mathbf s_n =  \mathbf x_i + r_{ijn}\mathbf s_{ijn} \in \Gamma_j^{\rm e}$, we can compute the following column vector $S$
\begin{equation}
\begin{array}{r@{}l}
\begin{aligned}
[S]_{ijn} & = \left( \widetilde{\mathcal S}_{\kappa,\Gamma_i}\chi_i^{\rm e}\partial_{\mathbf n}\psi_{0} \right)(\mathbf x_j + r_j\mathbf s_n) \\
&= \sum_{\ell_0 = 0}^{\ell_{\rm max}} \sum_{m=-\ell_0}^{\ell_0} c_{0,\ell_0 m}\, \left( \widetilde{\mathcal S}_{\kappa,\Gamma_i}Y_{\ell_0}^{m_0}\right)(\mathbf x_j + r_j\mathbf s_n)\\
& =  \sum_{\ell_0 = 0}^{\ell_{\rm max}} \sum_{m=-\ell_0}^{\ell_0} c_{0,\ell_0 m}\, \left( \frac{i'_{\ell_0}\left({ r_i}\right)}{i_{\ell_0}\left({}{r_i}\right)} -  \frac{k'_{\ell_0}\left({ r_i}\right)}{k_{\ell_0}\left({}{r_i}\right)} \right)^{-1} \,\frac{ k_{\ell_0}\left({ r_{ijn}}\right)}{k_{\ell_0}\left({}{r_i}\right)}\, Y_{\ell_0}^{m_0}(\mathbf s_{ijn}),
\end{aligned}
\end{array}
\end{equation} 
where
\begin{equation}
\begin{array}{r@{}l}
\begin{aligned}
c_{0,\ell_0 m_0}= \sum_{n = 1}^{N_{\rm leb}} w_n^{\rm leb} \chi_i^{\rm e}(\mathbf x_i + r_i\mathbf s_n) \partial_{\mathbf n}\psi_0(\mathbf x_i + r_i\mathbf s_n) Y_{\ell_0}^{m_0}(\mathbf s_n).
\end{aligned}
\end{array}
\end{equation}
This yields that
\begin{equation}\label{eq:F_0}
[F_0]_{j\ell m} = - \frac{\varepsilon_1}{\varepsilon_2} \left(\sum_{n = 1}^{N_{\rm leb}} w_n^{\rm leb} \chi_j^{\rm e}(\mathbf x_j +r_j\mathbf s_n) Y_{\ell}^m(\mathbf s_n) \sum_{i = 1}^M [S]_{ijn}\right).
\end{equation}

\end{document}